\documentclass[reqno,11pt]{amsart}
\usepackage{amscd,amssymb}

\theoremstyle{plain}
\newtheorem{Thm}[subsection]{Theorem}
\newtheorem{Cor}[subsection]{Corollary}
\newtheorem{Lem}[subsection]{Lemma}
\newtheorem{Prop}[subsection]{Proposition}
\newtheorem{Conj}[subsection]{Conjecture}

\theoremstyle{definition}
\newtheorem{Def}[subsection]{Definition}

\theoremstyle{remark}

\newtheorem{Rem}[subsection]{Remark}

\numberwithin{equation}{section}
\renewcommand{\rm}{\normalshape}

\newif\ifShowLabels
\ShowLabelstrue
\newdimen\theight
\def\TeXref#1{%
    \leavevmode\vadjust{\setbox0=\hbox{{\tt
        \quad\quad  {\small \rm #1}}}%
    \theight=\ht0
    \advance\theight by \lineskip
    \kern -\theight \vbox to
    \theight{\rightline{\rlap{\box0}}%
    \vss}%
    }}%

\ShowLabelsfalse% comment this out if labels should be printed

\renewcommand{\sec}[2]{\section{#2}\label{S:#1}%
    \ifShowLabels \TeXref{{S:#1}} \fi}
\newcommand{\ssec}[2]{\subsection{#2}\label{SS:#1}%
    \ifShowLabels \TeXref{{SS:#1}} \fi}

\newcommand{\refs}[1]{Section ~\ref{S:#1}}
\newcommand{\refss}[1]{Section ~\ref{SS:#1}}

\newcommand{\reft}[1]{Theorem ~\ref{T:#1}}

\newcommand{\refc}[1]{Corollary ~\ref{C:#1}}

\newcommand{\refe}[1]{\eqref{E:#1}}

\newenvironment{thm}[1]%
    { \begin{Thm} \label{T:#1}  \ifShowLabels \TeXref{T:#1} \fi }%
    { \end{Thm} }

\renewcommand{\th}[1]{\begin{thm}{#1} \sl }
\renewcommand{\eth}{\end{thm} }

\newenvironment{lemma}[1]%
    { \begin{Lem} \label{L:#1}  \ifShowLabels \TeXref{L:#1} \fi }%
    { \end{Lem} }
\newcommand{\lem}[1]{\begin{lemma}{#1} \sl}
\newcommand{\elem}{\end{lemma}}

\newenvironment{propos}[1]%
    { \begin{Prop} \label{P:#1}  \ifShowLabels \TeXref{P:#1} \fi }%
    { \end{Prop} }
\newcommand{\prop}[1]{\begin{propos}{#1}\sl }
\newcommand{\eprop}{\end{propos}}

\newenvironment{corol}[1]%
    { \begin{Cor} \label{C:#1}  \ifShowLabels \TeXref{C:#1} \fi }%
    { \end{Cor} }
\newcommand{\cor}[1]{\begin{corol}{#1} \sl }
\newcommand{\ecor}{\end{corol}}

\newenvironment{conjec}[1]%
    { \begin{Conj} \label{Co:#1}  \ifShowLabels \TeXref{Co:#1} \fi }%
    { \end{Conj} }
\newcommand{\conj}[1]{\begin{conjec}{#1} \sl }
\newcommand{\econj}{\end{conjec}}

\newenvironment{defeni}[1]%
    { \begin{Def} \label{D:#1}  \ifShowLabels \TeXref{D:#1} \fi }%
    { \end{Def} }
\newcommand{\defe}[1]{\begin{defeni}{#1} \sl }
\newcommand{\edefe}{\end{defeni}}

\newenvironment{remark}[1]%
    { \begin{Rem} \label{R:#1}  \ifShowLabels \TeXref{R:#1} \fi }%
    { \end{Rem} }
\newcommand{\rem}[1]{\begin{remark}{#1}}
\newcommand{\erem}{\end{remark}}

\newcommand{\eq}[1]%
    { \ifShowLabels \TeXref{E:#1} \fi
       \begin{equation} \label{E:#1} }
\newcommand{\eeq}{ \end{equation} }

\newcommand{\prf}{ \begin{proof} }
\newcommand{\epr}{ \end{proof} }

\newcommand\alp{\alpha}

\newcommand\del{\delta}     \newcommand\Del{\Delta}
\newcommand\eps{\varepsilon}

\newcommand\kap{\kappa}
\newcommand\lam{\lambda}        \newcommand\Lam{\Lambda}

\newcommand\ome{\omega}     

\newcommand\calA{{\mathcal{A}}}

\newcommand\calF{{\mathcal{F}}}
\newcommand\calG{{\mathcal{G}}}

\newcommand\calI{{\mathcal{I}}}

\newcommand\calK{{\mathcal{K}}}
\newcommand\calL{{\mathcal{L}}}
\newcommand\calM{{\mathcal{M}}}

\newcommand\calO{{\mathcal{O}}}

\newcommand\calU{{\mathcal{U}}}

\newcommand\calZ{{\mathcal{Z}}}

        \newcommand\bfC{{\mathbf C}}
        \newcommand\bfD{{\mathbf D}}

        \newcommand\bfS{{\mathbf S}}

        \newcommand\bfX{{\mathbf X}}

\newcommand\RR{\mathbb{R}}

\newcommand\PP{\mathbb{P}}
\renewcommand\AA{\mathbb{A}}
\renewcommand\SS{\mathbb{S}}
\newcommand\DD{\mathbb{D}}

\newcommand\GG{\mathbb{G}}

\newcommand\ZZ{\mathbb{Z}}

\newcommand\CC{\mathbb{C}}

 \newcommand\grg{{\mathfrak{g}}}
 \newcommand\grh{{\mathfrak{h}}}
 \newcommand\gri{{\mathfrak{i}}}
 \newcommand\grj{{\mathfrak{j}}}
 \newcommand\grk{{\mathfrak{k}}}
 \newcommand\grl{{\mathfrak{l}}}
 \newcommand\grm{{\mathfrak{m}}}
 \newcommand\grn{{\mathfrak{n}}}

 \newcommand\grq{{\mathfrak{q}}}

 \newcommand\grt{{\mathfrak{t}}}

\newcommand\sdp{\times \hskip -0.3em {\raise 0.3ex
\hbox{$\scriptscriptstyle |$}}} % semidirect product

\newcommand\ad{\operatorname{ad}}

\newcommand\Conv{\operatorname{Conv}}

\newcommand\End{\operatorname{End\,}}
\newcommand\Ext{\operatorname{Ext}}

\newcommand\GL{\operatorname{GL}}
\newcommand\Gr{\operatorname{Gr}}

\newcommand\Hom{\operatorname {Hom}}

\newcommand\id{\operatorname{id}}
\newcommand\Id{\operatorname{Id}}
\newcommand\im{\operatorname {im}}

\newcommand\Ker{\operatorname{Ker}}

\newcommand\pr{\operatorname{pr}}

\newcommand\Spec{\operatorname{Spec}}

\newcommand\Sym{\operatorname{Sym}}

\newcommand\uj{{\underline{j}}}

\newcommand\hatM{{\widehat{M}}}

\newcommand\x{\times}
\newcommand\ten{\otimes}

\newcommand{\ra}{\rangle}
\newcommand{\la}{\langle}

\newcommand{\Gra}{\operatorname{Graph}}
\newcommand\bt{\text{b}}

\newcommand\IC{\operatorname{IC}}
\newcommand\IH{\operatorname{IH}}

\newcommand{\nc}{\newcommand}
\nc{\renc}{\renewcommand}
\nc{\on}{\operatorname}

\nc\ol{\overline} \nc\wt{\widetilde} \nc\tboxtimes{\wt{\boxtimes}}
\renc{\SS}{{\mathbb S}} \renc{\DD}{{\mathbbD}}
\renewcommand{\AA}{{\mathbb A}}
\nc{\Fq}{{\mathbb F}_q} \nc{\Fqb}{\ol{{\mathbb F}_q}}
\nc{\Ql}{\ol{{\mathbb Q}_\ell}} \renc{\id}{\text{id}}
\nc\X{\mathcal X}

\renc{\Hom}{\on{Hom}} \nc{\Lie}{\on{Lie}} \nc{\Loc}{\on{Loc}}
\nc{\Pic}{\on{Pic}} \nc{\Bun}{\on{Bun}}
\nc{\Sh}{\on{Sh}}
\nc{\pos}{{\on{pos}}} \renc{\Conv}{\on{Conv}} \nc{\Sph}{\on{Sph}}
\renc{\Sym}{\on{Sym}}
\nc{\BunBb}{\overline{\Bun}_B} \nc{\Buno}{\overset{o}{\Bun}}
\nc{\BunPb}{{\overline{\Bun}_P}}
\nc{\BunBM}{\overline{\Bun}_{B(M)}}
\nc{\BunPbw}{{\widetilde{\Bun}_P}}
\nc{\BunBP}{\widetilde{\Bun}_{B,P}} \nc{\GUb}{\overline{G/U}}
\nc{\GUPb}{\overline{G/U(P)}}

\nc{\Hhom}{\underline{\on{Hom}}} \nc\syminfty{\on{Sym}^{\infty}}
\nc\lal{\ol{\lambda}} \nc\xl{\ol{x}} \nc\thl{\ol{\theta}}
\nc\nul{\ol{\nu}} \nc\mul{\ol{\mu}} \nc\Sum\Sigma
\nc{\hl}{\overset{\leftarrow}h}
\nc{\hr}{\overset{\rightarrow}h} \nc{\M}{{\mathcal M}}
\nc{\N}{{\mathcal N}} \nc{\F}{{\mathcal F}} \nc{\D}{{\mathcal D}}
\nc{\Q}{{\mathcal Q}} \nc{\Y}{{\mathcal Y}} \nc{\G}{{\mathcal G}}
\nc{\E}{{\mathcal E}} \nc{\CalC}{{\mathcal C}}
\nc\Dh{\widehat{\D}}

\nc{\C}{{\mathcal C}} \nc{\K}{{\mathcal K}}
\renewcommand{\H}{{\mathcal H}}

\nc{\T}{{\mathcal T}} \nc{\V}{{\mathcal V}} \renc{\P}{{\mathcal
P}} \nc{\A}{{\mathcal A}} \nc{\B}{{\mathcal B}} \nc{\U}{{\mathcal
U}}

\renc{\Gr}{\on{Gr}}

\nc{\frn}{{\check{\mathfrak u}(P)}} \nc{\p}{\mathfrak p}
\nc{\q}{\mathfrak q} \nc\f{{\mathfrak f}}

\nc{\qo}{{\mathfrak q}} \nc{\po}{{\mathfrak p}} \nc{\s}{{\mathfrak
s}} \nc\w{\text{w}}

\nc\mathi\iota \renc\Spec{\on{Spec}} \nc\Mod{\on{Mod}}
\nc{\tw}{\widetilde{\mathfrak t}} \nc{\pw}{\widetilde{\mathfrak
p}} \nc{\qw}{\widetilde{\mathfrak q}} \nc{\jw}{\widetilde j}

\nc{\I}{\mathcal I}

\nc{\lambdach}{{\check\lambda}} \nc{\Lambdach}{{\check\Lambda}{}}
\nc{\much}{{\check\mu}} \nc{\omegach}{{\check\omega}}
\nc{\nuch}{{\check\nu}} \nc{\etach}{{\check\eta}}
\nc{\alphach}{{\check\alpha}} \nc{\betach}{{\check\beta}}
\nc{\rhoch}{{\check\rho}} \nc{\ch}{{\check h}}

\nc{\Hb}{\overline{\H}}

\nc{\BA}{{\mathbb{A}}} \nc{\BC}{{\mathbb{C}}}
\nc{\BG}{{\mathbb{G}}} \nc{\BM}{{\mathbb{M}}}
\nc{\BN}{{\mathbb{N}}} \nc{\BP}{{\mathbb{P}}}
\nc{\BR}{{\mathbb{R}}} \nc{\BZ}{{\mathbb{Z}}}
\nc{\BS}{{\mathbb{S}}}

\nc{\CA}{{\mathcal{A}}} \nc{\CB}{{\mathcal{B}}}

\nc{\CE}{{\mathcal{E}}} \nc{\CF}{{\mathcal{F}}}
\nc{\CG}{{\mathcal{G}}} \nc{\CL}{{\mathcal{L}}}
\nc{\CM}{{\mathcal{M}}} \nc{\CN}{{\mathcal{N}}}
\nc{\CK}{{\mathcal{K}}} \nc{\CO}{{\mathcal{O}}}
\nc{\CP}{{\mathcal{P}}} \nc{\CQ}{{\mathcal{Q}}}
\nc{\CR}{{\mathcal{R}}} \nc{\CS}{{\mathcal{S}}}
\nc{\CT}{{\mathcal{T}}} \nc{\CU}{{\mathcal{U}}}
\nc{\CV}{{\mathcal{V}}} \nc{\CW}{{\mathcal{W}}}
\nc{\CZ}{{\mathcal{Z}}}

\nc{\cM}{{\check{\mathcal M}}{}} \nc{\csM}{{\check{\mathcal A}}{}}
\nc{\obM}{{\overset{\circ}{\mathbf M}}{}}
\nc{\oCA}{{\overset{\circ}{\mathcal A}}{}}
\nc{\obA}{{\overset{\circ}{\mathbf A}}{}}
\nc{\ooM}{{\overset{\circ}{M}}{}}
\nc{\osM}{{\overset{\circ}{\mathsf M}}{}}
\nc{\vM}{{\overset{\bullet}{\mathcal M}}{}}
\nc{\nM}{{\underset{\bullet}{\mathcal M}}{}}
\nc{\obD}{{\overset{\circ}{\mathbf D}}{}}
\nc{\cp}{{\overset{\circ}{\mathbf p}}{}}
\nc{\ofZ}{{\overset{\circ}{\mathfrak Z}}{}}
\nc{\fa}{{\mathfrak{a}}} \nc{\fb}{{\mathfrak{b}}}
\nc{\fg}{{\mathfrak{g}}} \nc{\fgl}{{\mathfrak{gl}}}
\nc{\fh}{{\mathfrak{h}}} \nc{\fj}{{\mathfrak{j}}}
\nc{\fm}{{\mathfrak{m}}} \nc{\fn}{{\mathfrak{n}}}
\nc{\fu}{{\mathfrak{u}}} \nc{\fp}{{\mathfrak{p}}}
\nc{\fr}{{\mathfrak{r}}} \nc{\fs}{{\mathfrak{s}}}
\nc{\fsl}{{\mathfrak{sl}}} \nc{\hsl}{{\widehat{\mathfrak{sl}}}}
\nc{\hgl}{{\widehat{\mathfrak{gl}}}}
\nc{\hg}{{\widehat{\mathfrak{g}}}}
\nc{\chg}{{\widehat{\mathfrak{g}}}{}^\vee}
\nc{\hn}{{\widehat{\mathfrak{n}}}}
\nc{\chn}{{\widehat{\mathfrak{n}}}{}^\vee}

\nc{\fA}{{\mathfrak{A}}} \nc{\fB}{{\mathfrak{B}}}
\nc{\fD}{{\mathfrak{D}}} \nc{\fE}{{\mathfrak{E}}}
\nc{\fF}{{\mathfrak{F}}} \nc{\fG}{{\mathfrak{G}}}
\nc{\fK}{{\mathfrak{K}}} \nc{\fL}{{\mathfrak{L}}}
\nc{\fM}{{\mathfrak{M}}} \nc{\fN}{{\mathfrak{N}}}
\nc{\fP}{{\mathfrak{P}}} \nc{\fU}{{\mathfrak{U}}}
\nc{\fV}{{\mathfrak{V}}} \nc{\fZ}{{\mathfrak{Z}}}

\nc{\bb}{{\mathbf{b}}} \nc{\bc}{{\mathbf{c}}}
\nc{\bd}{{\mathbf{d}}} \nc{\be}{{\mathbf{e}}}
\nc{\bj}{{\mathbf{j}}} \nc{\bn}{{\mathbf{n}}}
\nc{\bp}{{\mathbf{p}}} \nc{\bq}{{\mathbf{q}}}
\nc{\bu}{{\mathbf{u}}} \nc{\bv}{{\mathbf{v}}}
\nc{\bx}{{\mathbf{x}}} \nc{\bs}{{\mathbf{s}}}
\nc{\by}{{\mathbf{y}}} \nc{\bw}{{\mathbf{w}}}
\nc{\bA}{{\mathbf{A}}} \nc{\bK}{{\mathbf{K}}}
\nc{\bB}{{\mathbf{B}}} \nc{\bC}{{\mathbf{C}}}
\nc{\bD}{{\mathbf{D}}} \nc{\bH}{{\mathbf{H}}}
\nc{\bM}{{\mathbf{M}}} \nc{\bN}{{\mathbf{N}}}
\nc{\bV}{{\mathbf{V}}} \nc{\bW}{{\mathbf{W}}}
\nc{\bX}{{\mathbf{X}}} \nc{\bZ}{{\mathbf{Z}}}
\nc{\bS}{{\mathbf{S}}}

\nc{\sA}{{\mathsf{A}}} \nc{\sB}{{\mathsf{B}}}
\nc{\sC}{{\mathsf{C}}} \nc{\sD}{{\mathsf{D}}}
\nc{\sF}{{\mathsf{F}}} \nc{\sK}{{\mathsf{K}}}
\nc{\sM}{{\mathsf{M}}} \nc{\sO}{{\mathsf{O}}}
\nc{\sQ}{{\mathsf{Q}}} \nc{\sP}{{\mathsf{P}}}
\nc{\sZ}{{\mathsf{Z}}} \nc{\sfp}{{\mathsf{p}}}
\nc{\sr}{{\mathsf{r}}} \nc{\sg}{{\mathsf{g}}}
\nc{\sff}{{\mathsf{f}}} \nc{\sfb}{{\mathsf{b}}}
\nc{\sfc}{{\mathsf{c}}} \nc{\sd}{{\mathsf{d}}}

\nc{\BK}{{\bar{K}}}

\nc{\tA}{{\widetilde{\mathbf{A}}}}
\nc{\tB}{{\widetilde{\mathcal{B}}}}
\nc{\tg}{{\widetilde{\mathfrak{g}}}} \nc{\tG}{{\widetilde{G}}}
\nc{\TM}{{\widetilde{\mathbb{M}}}{}}
\nc{\tO}{{\widetilde{\mathsf{O}}}{}}
\nc{\tU}{{\widetilde{\mathfrak{U}}}{}} \nc{\TZ}{{\tilde{Z}}}
\nc{\tx}{{\tilde{x}}} \nc{\tbv}{{\tilde{\bv}}}
\nc{\tfP}{{\widetilde{\mathfrak{P}}}{}} \nc{\tz}{{\tilde{\zeta}}}
\nc{\tmu}{{\tilde{\mu}}}

\nc{\urho}{\underline{\rho}} \nc{\uB}{\underline{B}}
\nc{\uC}{{\underline{\mathbb{C}}}} \nc{\ui}{\underline{i}}
\renc{\uj}{\underline{j}} \nc{\ofP}{{\overline{\mathfrak{P}}}}
\renc{\eps}{\varepsilon} \nc{\hrho}{{\hat{\rho}}}

\nc{\one}{{\mathbf{1}}} \nc{\two}{{\mathbf{t}}}

\nc{\Rep}{{\mathop{\operatorname{\rm Rep}}}}
\nc{\Tot}{{\mathop{\operatorname{\rm Tot}}}}
\renc{\Ker}{{\mathop{\operatorname{\rm Ker}}}}
\nc{\Hilb}{{\mathop{\operatorname{\rm Hilb}}}}
\renc{\End}{{\mathop{\operatorname{\rm End}}}}
\renc{\Ext}{{\mathop{\operatorname{\rm Ext}}}}
\nc{\CHom}{{\mathop{\operatorname{{\mathcal{H}}\it om}}}}
\renc{\GL}{{\mathop{\operatorname{\rm GL}}}}
\nc{\gr}{{\mathop{\operatorname{\rm gr}}}}
\renc{\Id}{{\mathop{\operatorname{\rm Id}}}}
\nc{\de}{{\mathop{\operatorname{\rm def}}}}
\nc{\length}{{\mathop{\operatorname{\rm length}}}}
\nc{\Cliff}{{\mathsf{Cliff}}}
\nc{\Fl}{\on{Fl}} \nc{\Fib}{{\mathsf{Fib}}}
\nc{\Coh}{{\mathsf{Coh}}} \nc{\FCoh}{{\mathsf{FCoh}}}

\nc{\reg}{{\text{\rm reg}}}

\nc{\cplus}{{\mathbf{C}_+}} \nc{\cminus}{{\mathbf{C}_-}}
\nc{\cthree}{{\mathbf{C}_*}} \nc{\Qbar}{{\bar{Q}}}

\nc{\bh}{{\bar{h}}} \nc{\bOmega}{{\overline{\Omega}}}

\nc{\seq}[1]{\stackrel{#1}{\sim}} \nc\QM{{\mathcal {QM}}}

\nc{\chH}{\check H} \nc{\chM}{\check M} \nc{\aff}{{\on{aff}}}

\nc{\chh}{\check \grh}

\renewcommand\chn{\check \grn}

\newcommand{\ocM}{{\overline\calM}}

\renewcommand\chg{\check \grg}
\newcommand\ev{\operatorname{ev}}

\newcommand\chd{\check d}
\newcommand\chK{\check K}
\nc\chT{\check T}
\renewcommand\Ext{\operatorname{Ext}}

\begin{document}
\title[Instanton counting via affine Lie algebras]{Instanton counting via affine Lie algebras I:
equivariant J-functions of (affine) flag manifolds and Whittaker vectors}
\author{Alexander Braverman}
\email{braval@math.harvard.edu}
\address{Department of Mathematics, Harvard University, 1 Oxford
street, Cambridge MA 02138 USA}
\thanks{This research has been partially supported by the NSF}

\begin{abstract}Let $\grg$ be a simple complex Lie algebra,
$G$ - the corresponding simply connected group; let also
$\grg_{\aff}$ be the corresponding untwisted affine Lie algebra.
For a parabolic subgroup $P\subset G$ we introduce a generating
function $\calZ_{G,P}^{\aff}$ which roughly speaking counts framed
$G$-bundles on $\PP^2$  endowed with a $P$-structure on the
horizontal line (the formal definition uses the corresponding
Uhlenbeck type compactifications studied in \cite{bfg}). In the
case $P=G$ the function $\calZ_{G,P}^{\aff}$ coincides with
Nekrasov's partition function introduced in \cite{nek} and studied
thoroughly in \cite{neok} and \cite{nayo} for $G=SL(n)$. In the
"opposite case" when $P$ is a Borel subgroup of $G$ we show that
$\calZ_{G,P}^{\aff}$ is equal (roughly speaking) to the Whittaker
matrix coefficient in the universal Verma module for the Lie
algebra $\chg_{\aff}$ -- the Langlands dual Lie algebra of
$\grg_{\aff}$. This clarifies somewhat the connection between
certain asymptotic of $\calZ_{G,P}^{\aff}$ (studied in {\it loc.
cit.} for $P=G$) and the classical affine Toda system. We also
explain why the above result gives rise to a calculation of (not
yet rigorously defined) equivariant quantum cohomology ring of the
affine flag manifold associated with $G$. In particular, we
reprove the results of \cite{giki} and \cite{kim} about quantum
cohomology of ordinary flag manifolds using methods which are
totally different from {\it loc. cit.}

We shall show in a subsequent publication how this allows one to connect certain asymptotic of the
function $\calZ_{G,P}^{\aff}$ with {\it the Seiberg-Witten prepotential} (cf. \cite{bret},
thus proving the main conjecture of
\cite{nek} for an arbitrary gauge group $G$ (for $G=SL(n)$ it has been proved in \cite{neok} and \cite{nayo}
by other methods.
\end{abstract}
\maketitle
%--------------------------------------------------------------------

\sec{}{Introduction}

\ssec{}{The partition function}
This paper has grown out of a
(still unsuccessful) attempt to understand the following object.
Let $K$ be a simply connected connected compact Lie group and let
$d$ be a non-negative integer. Denote by $\calM_K^d$ the moduli
space of (framed) $K$-instantons on $\RR^4$  of second Chern class
$-d$. This space can naturally be embedded into a larger {\it
Uhlenbeck space} $\calU_K^d$. Both spaces admit a natural action of
the group $K$ (by changing the framing at $\infty$) and the torus
$(S^1)^2$ acting on $\RR^4$ after choosing an identification
$\RR^4\simeq \CC^2$. Moreover the maximal torus of $K\x (S^1)^2$
has unique fixed point on $\calU_K^d$. Thus we may consider (cf.
\cite{nek}, \cite{nayo} and \refs{equivariant integration} for
precise definitions) the {\it equivariant integral}
$$
\int\limits_{\calU^d_K}1^d
$$
of the unit $K\times (S^1)^2$-equivariant cohomology class (which we
denote by $1^d$) over
$\calU_K^d$; the integral takes values in the field $\calK$ which is
the field of fractions of the algebra $\calA=H^*_{K\x
(S^1)^2}(pt)$ \footnote{In this paper we always consider
cohomology with complex coefficients}. Note that $\calA$ is
canonically isomorphic to the algebra of polynomial functions on
$\grk\x\RR^2$ (here $\grk$ denotes the Lie algebra of $K$) which
are invariant with respect to the adjoint action of $K$ on $\grk$.
Thus each $\int\limits_{\calU^d_K}1^d$ may naturally be regarded as a
rational function of $a\in\grk$ and $(\eps_1,\eps_2)\in \RR^2$.

Consider now the generating function
$$
\calZ=\sum\limits_{d=0}^\infty Q^d \int\limits_{\calU_K^d} 1^d.
$$
It can (and should) be thought of as a function of the variables
$\grq$ and $a,\eps_1,\eps_2$ as before. In \cite{nek} it was
conjectured that the first term of the asymptotic in the limit
$\underset{\eps_1,\eps_2\to 0}\lim \ln \calZ$ is closely related to
{\it Seiberg-Witten prepotential} of $K$. For $K=SU(n)$ this
conjecture has been proved in \cite{neok} and \cite{nayo}. Also in
\cite{nek} an explicit combinatorial expression for $\calZ$ has
been found.
%------------------------------------------------------------------------
\ssec{}{Algebraic version}In this paper we want to generalize the
definition of the function $\calZ$ in several directions and
compute it in some of these new cases. First of all, it will be
convenient for us to make the whole situation completely
algebraic.

Namely, let $G$ be a complex semi-simple algebraic group whose
maximal compact subgroup is isomorphic to $K$. We shall denote by
$\grg$ its Lie algebra. Let also $\bfS=\PP^2$ and denote by
$\bfD_\infty\subset \bfS$ the "straight line at $\infty$"; thus
$\bfS\backslash\bfD_\infty=\AA^2$. It is well-known that $\calM^d_K$
is isomorphic to the moduli space $\Bun_G^d(\bfS,\bfD_\infty)$ of
principal $G$-bundles on $\bfS$ endowed with a trivialization on
$\bfD_\infty$ of second Chern class $-d$. When it does not lead to
a confusion we shall write $\Bun_G$ instead of
$\Bun_G(\bfS,\bfD_\infty)$. The algebraic analog of $\calU^d_K$ has
been constructed in \cite{bfg}; we denote this algebraic variety
by $\calU_G^d$. This variety is endowed with a natural action on
$G\x (\CC^*)^2$.
%---------------------------------------------------------------------------------------------------------------
\ssec{}{Parabolic generalization of the partition function}Let
$\bfC\subset\bfS$ denote the "horizontal line". Choose a parabolic
subgroup $P\subset G$. Let $\Bun_{G,P}$ denote the moduli space of
the following objects:

1) A principal $G$-bundle $\calF_G$ on $\bfS$;

2) A trivialization of $\calF_G$ on $\bfD_\infty\subset \bfS$;

3) A reduction of $\calF_G$ to $P$ on $\bfC$ compatible with the
trivialization of $\calF_G$ on $\bfC$.

Let us describe the connected components of $\Bun_{G,P}$. Let $M$
be the Levi group of $P$. Denote by $\chM$ the {\it Langlands
dual} group of $M$ and let $Z(\chM)$ be its center. We denote by
$\Lam_{G,P}$ the lattice of characters of $Z(\chM)$. Let also
$\Lam_{G,P}^{\aff}=\Lam_{G,P}\x \ZZ$ be the lattice of characters
of $Z(\chM)\x\CC^*$. Note that $\Lam_{G,G}^{\aff}=\ZZ$.

The lattice $\Lam^{\aff}_{G,P}$ contains canonical semi-group
$\Lam^{\aff,\pos}_{G,P}$ of positive elements (cf. \cite{bfg} and
\refs{main}). It is not difficult to see that the connected
components of $\Bun_{G,P}$ are parameterized  by the elements of
$\Lam_{G,P}^{\aff,\pos}$:
$$
    \Bun_{G,P}=\bigcup\limits_{\theta_{\aff}\in\Lam_{G,P}^{\aff,\pos}}
    \Bun_{G,P}^{\theta_{\aff}}.
$$

Typically, for $\theta_{\aff}\in \Lam_{G,P}^{\aff}$ we shall write
$\theta_\aff=(\theta,d)$ where $\theta\in \Lam_{G,P}$ and $d\in \ZZ$.

Each $\Bun_{G,P}^{\theta_{\aff}}$ is naturally acted on by $P\x(\CC^*)^2$;
by embedding $M$ into $P$ we get an action of $M\x(\CC^*)^2$ on
$\Bun_{G,P}^{\theta_{\aff}}$. In \cite{bfg} we define for each
$\theta_{\aff}\in\Lam_{G,P}^{\aff,\pos}$ certain Uhlenbeck scheme
$\calU_{G,P}^{\theta_{\aff}}$ which contains $\Bun_{G,P}^{\theta_{\aff}}$ as a dense
open subset. The scheme $\calU_{G,P}^{\theta_{\aff}}$ still admits an
action of $M\x(\CC^*)^2$.

We want to do some equivariant intersection theory on the spaces
$\calU_{G,P}^{\theta_{\aff}}$. For this let us denote by $\calA_{M\x(\CC^*)^2}$
the algebra $H^*_{M\x (\CC^*)^2}(pt,\CC)$. Of course this is just
the algebra of $M$-invariant polynomials on $\grm\x\CC^2$.
Let also $\calK_{M\x(\CC^*)^2}$ be its field of fractions. We can think
about elements of $\calK_{M\x (\CC^*)^2}$ as rational functions on
$\grm\x\CC^2$ which are invariant with respect to the adjoint action.

Let $T\subset M$ be a maximal torus. Then one can show that
$(\calU_{G,P}^{\theta_{\aff}})^{T\x(\CC^*)^2}$ consists of one point.
This guarantees that we may consider the integral
$\int\limits_{\calU_{G,P}^{\theta_{\aff}}}1_{G,P}^{\theta_{\aff}}$ where
$1_{G,P}^{\theta_{\aff}}$
denotes the unit class in
$H^*_{M\x(\CC^*)^2}(\calU_{G,P}^{\theta_{\aff}},\CC)$.
The result can be thought
of as a rational function on $\grm\x\CC^2$ which is invariant with
respect to the adjoint action of $M$. Define
\eq{partition-affine}
    \calZ_{G,P}^{\aff}=\sum\limits_{\theta\in\Lam_{G,P}^{\aff}} \grq_{\aff}^{\theta_{\aff}}
    \ \int\limits_{\calU_{G,P}^{\theta_{\aff}}}1_{G,P}^{\theta_{\aff}}.
\end{equation}
One should think of $\calZ_{G,P}^{\aff}$ as a formal power series
in $\grq_{\aff}\in Z(\chM)\x\CC^*$ with values in the space of ad-invariant
rational functions on $\grm\x\CC^2$. Typically, we shall write
$\grq_{\aff}=(\grq,Q)$ where $\grq\in Z(\chM)$ and $Q\in\CC^*$.
Also we shall denote an element of $\grm\x\CC^2$ by
$(a,\eps_1,\eps_2)$ or (sometimes it will be more convenient) by
$(a,\hbar,\eps)$ (note that for general $P$ (unlike in the case $P=G$)
the function $\calZ_{G,P}^{\aff}$ is not symmetric with respect to
switching
$\eps_1$ and $\eps_2$).

%-----------------------------------------------------------
\ssec{}{Interpretation via maps and the "finite-dimensional" analog}
Choose now another smooth
projective curve $\bfX$ of genus $0$ and with two marked points
$0_\bfX,\infty_\bfX$. Choose also a coordinate $x$ on $\bfX$ such
that $x(0_\bfX)=0$ and $x(\infty_\bfX)=0$. Let us denote by
$\CG_{G,P,\bX}$ the scheme classifying triples
$(\F_G,\beta,\gamma)$, where

1) $\F_G$ is a principal $G$-bundle on $\bfX$;

2) $\beta$ is a trivialization of $\F_G$ on the formal
neighborhood of $\infty_\bfX$;

3)  $\gamma$ is  a reduction to $P$ of the fiber of $\F_G$ at
$0_\bX$.

We shall usually omit $\bfX$ from the
notations.
We shall also write $\calG_G^{\aff}$ for
$\calG_{G,G}^{\aff}$.

Let $e_{G,P}^{\aff}\in\calG_{G,P}^{\aff}$ denote the point
corresponding to the trivial $\calF_G$ with the natural $\beta$
and $\gamma$. It is explained in \cite{bfg} that the variety
$\Bun_{G,P}$ is canonically isomorphic to the scheme classifying
{\it based maps} from $(\bfC,\infty_\bfC)$ to
$(\calG_{G,P}^{\aff},e_{G,P}^{\aff})$ (i.e. maps from $\bfC$
to $\calG_{G,P}$ sending $\infty_\bfC$ to
$e_{G,P}^{\aff}$).

The scheme $\calG_{G,P}^{\aff}$ may (and should) be thought of as
a partial flag variety for $\grg_{\aff}$. Thus we may consider the
following "finite-dimensional" analog of the above problem: for
$G$ and $P$ as above let $\calG_{G,P}=G/P$. Let also
$e_{{G,P}}\in\calG_{G,P}$ denote the image of $e\in G$.
Clearly $e_{G,P}$ is stable under the action of $P$ on
$\calG_{G,P}$. Let $^{\bt}\calM_{G,P}$ denote the moduli space of based
maps from $(\bfC,\infty_{\bfC})$ to
$(\calG_{G,P},e_{{G,P}})$, i.e. the moduli space of maps
$\bfC\to\calG_{G,P}$ which send $\infty_{\bfC}$ to $e_{G,P}$.
This space is acted on by the group $M\x\CC^*$.
Also $^{\bt}\calM_{G,P}$ is a union of connected components
$^{\bt}\calM_{G,P}^{\theta}$ where $\theta$ lies in a certain
sub-semi-group
$\Lam_{G,P}^{\pos}$ of $\Lam_{G,P}$. For each
$\theta$ as above one can also consider the space of {\it based
quasi-maps} (or {\it Zastava space} in the terminology of
\cite{bgfm}, \cite{fm} and \cite{ffkm}) which we denote by
$^{\bt}\QM^{\theta}_{G,P}$. This is an affine algebraic variety
containing $^{\bt}\calM^{\theta}_{G,P}$ as a dense open subset.
We also denote by $\calM_{G,P}^{\theta},\QM_{G,P}^{\theta}$ the
corresponding
spaces of all (unbased) quasi-maps.

We now introduce the ``finite-dimensional'' analog of the partition
function
\refe{partition-affine}. As before let
$$
\calA_{M\x\CC^*}=H^*_{M\x\CC^*}(pt,\CC)
$$
and  denote by $\calK_{M\x\CC^*}$ its field of fractions. Let also
$1_{G,P}^{\theta}$ denote the
unit class in the
$M\x \CC^*$-equivariant cohomology of $^{\bt}\QM_{G,P}^{\theta}$.
Then we define
\eq{partition-finite}
       \calZ_{G,P}=\sum\limits_{\theta\in\Lam_{G,P}^{\theta}}
\grq^{\theta} \ \int\limits_{^{\bt}\QM_{G,P}^{\theta}}1_{G,P}^{\theta}.
\end{equation}
This is a formal series in $q\in Z(\chM)$ with values in the field
$\calK_{M\x\CC^*}$ of $M$-invariant rational functions on $\grm\x\CC$.

In fact the function $\calZ_{G,P}$ is a familiar object in Gromov-Witten theory: we
shall
see later (\reft{j versus z}) that up to a simple factor $\calZ_{G,P}$
is the so called {\it equivariant $J$-function} of $\calG_{G,P}$
(cf. \refss{j-function} for the details).
%---------------------------------------------------------
\ssec{borel-int}{The Borel case} We believe that it should be possible to
express the function $\calZ_{G,P}$ (resp. the function
$\calZ_{G,P}^{\aff}$) in terms of representation theory of the Lie
algebra $\chg$ (resp. $\chg_{\aff}$) -- by the definition this is a
Lie algebra whose root system is dual to that of $\grg$ (resp. to that
of $\grg_{\aff}$). One of the motivations for
this comes from the main results of \cite{bgfm} and \cite{bfg}
where such a description is found for the intersection cohomology
of the varieties $^{\bt}\QM_{G,P}^{\theta}$ and
$\calU_{G,P}^{\theta_{\aff}}$ (in this paper we are going to adopt
an intersection cohomology approach to the partition functions
$\calZ_{G,P}$ and $\calZ_{G,P}^{\aff}$; this is explained
carefully in \refs{equivariant integration}). One of the main
results of this paper gives such a calculation of the functions
$\calZ_{G,B}$ and $\calZ_{G,B}^{\aff}$ where $B\subset G$ is a
Borel subgroup of $G$. Roughly speaking we show that $\calZ_{G,B}$
(resp. $\calZ_{G,B}^{\aff}$) is equal to {\it Whittaker matrix
coefficient} of the Verma module over $\chg$ (resp. over
$\chg_{\aff}$) whose lowest weight given by $\frac{a}{\hbar}+\rho$
(resp. $\frac{(a,\eps_1)}{\eps_2}+\rho_{\aff}$ where $a,\hbar,\eps_1$ and $\eps_2$
are as in the previous subsection (here we regard $(a,\eps_1)$ as a
weight
for the dual affine algebra $\chg_{\aff}$; this is explained carefully
in Section 3).
The precise formulation of this result is given by
\reft{main} and
\refc{main-c}.

As a corollary we get that the function
$\grq^{\frac{a}{\hbar}}\calZ_{G,B}$ is an eigen-function of
the{\it quantum Toda hamiltonians}
 associated with $\chg$
with eigen-values determined (in the natural way) by $a$
 (we refer the reader to
\cite{etingof} for the definition of (affine) Toda integrable system
and its relation with Whittaker functions).
In the affine case one can also show that
$\grq^{\frac{a}{\hbar}}\calZ_{G,B}^{\aff}$
is an eigen-function of a certain differential operator which has
order 2 (``non-stationary analog'' of the affine quadratic Toda hamiltonian).
We shall show how this allows to compute the asymptotics of {\it all}
the functions $\calZ_{G,P}^{\aff}$ when $\eps_1,\eps_2\to 0$ in
another publication (cf. \cite{bret}).

Let us explain the idea of the proof of the above statement (we shall do it
in the ``finite-dimensional'' (i.e. non-affine) case; in the affine
case the proof is similar). First for a scheme $Y$ let $\IH^*(Y)$
denote the intersection cohomology of $Y$ with complex coefficients
(similarly if a group $L$ acts on $Y$ we denote by $\IH^*_L(Y)$ the
corresponding $L$-equivariant  intersection cohomology). Let now
$$
\IH_{G,P}^{\theta}=\IH^*_{M\x(\CC^*)^2}
(^{\bt}\QM_{G,P}^{\theta})\underset{\calA_{M\x\CC^*}}
\ten \calK_{M\x\CC^*}.
$$
(it is clear that $\IH^*_{M\x(\CC^*)^2}
(^{\bt}\QM_{G,P}^{\theta})$ has a natural $\calA_{M\x\CC^*}$-module structure).

Let also
$$
\IH_{G,P}
=\bigoplus\limits_{\theta\in \Lam_{G,P}^{\pos}}\IH_{G,P}^{\theta}
$$
Each $\IH_{G,P}^{\theta}$ is a finite-dimensional vector space over
$\calK_{M\x\CC^*}$
endowed with a (non-degenerate) Poincar\'e pairing $\la\cdot,\cdot\ra_{G,P}^{\theta}$
taking values in
$\calK_{M\x\CC^*}$.
.

Let us now specialize to the case $P=B$.
In \refs{action} we show that the Lie algebra $\chg$
acts naturally on $\IH_{G,B}$. Moreover, this action has the following
properties. First of all, let us denote by
$\la\cdot,\cdot\ra_{G,B}$ the direct sum of the pairings
$(-1)^{\la \theta,\check{\rho}\ra}\la\cdot,\cdot\ra_{G,B}^{\theta}$.
\footnote{here $\check{\rho}$ denotes the half-sum of the positive roots of $\grg$}

Recall that the Lie algebra $\chg$ has its triangular
decomposition
$\chg=\chn_+\oplus\chh\oplus\chn_-$. Let $\kap:\chg\to\chg$ denote the
Cartan anti-involution which interchanges $\chn_+$ and $\chn_-$ and
acts as identity on $\chh$.
For each $\lam\in\grh=(\chh)^*$ we denote by $M(\lam)$ the
corresponding
Verma module with lowest weight $\lam$; this is a module generated by
a vector $v_\lam$ with (the only) relations
$$
t(v_\lam)=\lam(t)v_\lam\quad\text{for $t\in\chh$ and}\quad
n(v_\lam)=0\quad\text{for $n\in\chn_-$}.
$$

Then:

1) $\IH_{G,B}$ (with the
above action) becomes isomorphic to $M(\lam)$.

2) $\IH_{G,B}^{\theta}\subset \IH_{G,B}$ is the
$\frac{a}{\hbar}+\rho+\theta$-weight space of $\IH_{G,B}$.

3) For each $g\in\chg$ and $v,w\in \IH_{G,B}$ we have
$$
\la g(v),w\ra_{G,B}=\la v,\kap(g)w\ra_{G,B}.
$$

4) The  vector
$\sum_{\theta}1_{G,B}^{\theta}$ (lying is some completion of
$\IH_{G,B}$) is a Whittaker vector (i.e. a $\grn_-$-eigen-vector) for the above action.

It is easy to see that the assertions 1-4 imply the above
representation-theoretic description of $\calZ_{G,B}$ (for this one
has
to adopt an ``intersection cohomology'' point of view at the above
equivariant integrals; this is explained in \refs{equivariant
integration}).

In fact, certain parts of the above statement have been known
before: for example in \cite{ffkm} the authors compute the
dimensions of the spaces $\IH^{\theta}_{G,B}$ and the answer
agrees with 1 and 2 above (in fact that these dimensions are known
even if we replace $^{\bt}\QM_{G,B}^{\theta}$ by any
$^{\bt}\QM_{G,P}^{\theta}$ - cf. \cite{bgfm}). Similarly, the
computation of dimensions of $\IH_{G,B}^{\theta_{\aff}}$ (defined
analogously to $\IH_{G,B}^{\theta}$) follows from \cite{bfg}. Also
the construction of the $\chg$-action on $\IH_{G,B}$ is very close
to the construction in Section 4 of \cite{ffkm}.
%-----------------------------------------------------------------
\ssec{j-function}{Connection with $J$-functions and quantum cohomology of
flag manifolds}
Let us explain how the above result is connected
with the known results on quantum cohomology of flag manifolds.

First, let us recall the general set-up (the main reference is \cite{giv}).
Let $X$ be any smooth
projective variety over $\CC$. Assume that $X$ is homogeneous,
i.e. that the tangent sheaf of $X$ is generated by its global
sections. Then for any $\beta\in H_2(X,\ZZ)$ one may consider the
moduli space $\ocM_{0,n}(X,\beta)$ of stable maps from a curve
$\Sigma$ of genus zero with $n$ marked points to $X$ of degree
$\beta$.  For every $\beta$ as above let $\ev_{\beta}$ denote the
evaluation map map at the marked point from $\ocM_{0,1}(X,\beta)$
to $X$.

The space $\ocM_{0,1}(X,\beta)$ admits canonical line bundle
$\calL_\beta$ whose fiber over every stable map $C\to X$ is equal
to the cotangent space to $C$ at the marked point. We denote by
$c_{\beta}$ its first Chern class.

In the above set-up one defines the $J$-function $J^X$ of $X$ in
the following way:
\eq{j-function} %
J^X(t,\hbar)=e^{\frac{t}{\hbar}}\sum\limits_{\beta\in
H_2(X,\ZZ)}e^{\la
t,\beta\ra}(\ev_\beta)_*(\frac{1}{\hbar(c_{\beta}+\hbar)})
\end{equation}
where $t\in H^2(X,\CC)$\footnote{For general $X$ the definition of the
push-forward is somewhat tricky - it uses the so called virtual fundamental class.
However, for $X$ of the form
${\mathcal G}_{G,P}$ this difficulty is not present.}. The function $J^X$ takes values in
$H^*(X,\CC)\ten \CC((\hbar^{-1}))$ which we may think of as a
localization of a
completion of
of $H_{\CC^*}(X,\CC)$ ($\CC^*$ acts trivially on $X$) over
$\CC[\hbar]=H^*_{\CC^*}(pt)$. It is explained in
\cite{giv} that by looking at the differential equations satisfied
by the function $J^X$ one may compute explicitly the small quantum
cohomology ring of $X$. Putting $q=e^t$ one may (formally)
write
$$
  J^X=
  \grq^{\frac{1}{\hbar}}\sum\limits_{\beta\in
  H_2(X,\ZZ)}\grq^{\beta}(\ev_\beta)_*(\frac{1}{\hbar(c_{\beta}+\hbar)})
$$
In the case when $X$ is acted on by a reductive group $G$ we may
consider the equivariant analog of $J^X$ which we shall denote by
$J^X_G$. In this case the function $J^X_G$ takes values in
$H^*_{G\x\CC^*}(X,\CC)\underset{\CC[\hbar]}\ten \CC((\hbar^{-1}))$. By looking at the differential equations
satisfied by $J^X_G$ one may compute the equivariant small quantum
cohomology ring of $X$.

Assume now that $X=\calG_{G,P}$. Note that in this case we have
the natural identification
$$
H^*_{G\x \CC^*}(\calG_{G,P},\CC)=H^*_{M\x\CC^*}(pt)
$$
Note also that $H_2(\calG_{G,P},\ZZ)=\Lam_{G,P}$ and
$H^2(\calG_{G,P})=\Lie Z(\chM)$.

Taking all these identifications into account we now claim the
following
%-------------------------------------------------------------------
\th{j versus z}For any $G$ and $P$ as above we have
$$
J_G^{\calG_{G,P}}=\grq^{\frac{a}{\hbar}}\calZ_{G,P}.
$$
\eth
%----------------------------------------------------------------
\reft{j versus z} is proved in \refs{quantum cohomology}.

Note that in view of the results announced in the previous section
it follows that $J_{\calG_{G,B}}^G$ can be interpreted via Whittaker
matrix coefficients of Verma modules and in particular it is an
eigen-function
of the Toda hamiltonians for $\chg$. The latter result is due to
B.~Kim (cf. \cite{kim}); in fact (the ``non-affine'' part of) this
paper
may be viewed as a conceptual explanation of Kim's result.

\reft{j versus z} also allows us to think about
$\grq^{\frac{a}{\hbar}}\calZ_{G,P}^{\aff}$
as the equivariant (with respect to the loop group) $J$-function of
$\calG_{G,P}^{\aff}$. The theory of Gromov-Witten invariants of
$\calG_{G,P}^{\aff}$ is  not yet constructerd. However we see from the above
that we have a natural candidate for the notion of $J$-function of
$\calG_{G,P}$. In particular, the main result of this paper gives a
representation-theoretic interpretation of this $J$-function. This may
be viewed as an explanation of the results of \cite{go} and
\cite{mare}
(where the authors  compute the quantum cohomology of
$\calG_{G,P}$ making some assumptions about its existence and properties).
%---------------------------------------------------------------
%--------------------------------------------------
\ssec{}{Organization of the paper}This paper is organized as follows.
In \refs{equivariant integration} we fix the notation and recall the
basic facts about equivariant integration. In \refs{main} we give the
precise formulation of the main result of this paper discussed above in
\refss{borel-int} (\reft{main}). \refs{sl2} is devoted to an explicit
proof of the main result for $G=SL(2)$. In \refs{action} we prove
\reft{main}. Finally, \refs{quantum cohomology} is devoted to the proof of \reft{j
versus z}.
%--------------------------------------------------------------------
\ssec{}{Acknowledgments}We are grateful to N.~Nekrasov for his patient
explanation of the contents of \cite{nek} and \cite{neok} and numerous
discussions on the subjects. We also thank T.~Coates, P.~Etingof, M.~Finkelberg, D.~Gaitsgory,
D.~Kazhdan, H.~Nakajima and A.~Okounkov
for interesting conversations and help at various stages of this work.
%---------------------------------------------------------------------
\sec{equivariant integration}{Equivariant integration}
%-----------------------------------------------------------
\ssec{}{}Let $L$ be a reductive algebraic group over $\CC$. Set
${\calA_L}=H^*_L(pt,\CC)$. It is well known that ${\calA_L}$ is
canonically isomorphic to the algebra of invariant polynomials on
the Lie algebra $\grl$ of $L$. Let also ${\calK_L}$ denote the
field of fractions of ${\calA_L}$. One can identify ${\calK_L}$
with the field of invariant rational functions on $\grl$. We shall
say that an element $f\in{\calK_L}$ is homogeneous of degree
$d\in\ZZ$ if $f$ is a homogeneous rational function of degree $d$
on $\grl$.

\ssec{}{Equivariant (co)homology}
Let $Y$ be a scheme of finite type over $\CC$ endowed with an
action of $L$. We denote by $\ome_Y$ the dualizing complex of
$Y$ (in the category of constructible sheaves). We will be interested
in the equivariant cohomology $H^*_L(Y,\ome_Y)$ (it is sometimes called
the {\it equivariant Borel-Moore homology} of $Y$).
For any closed $L$-invariant subscheme $Z\subset Y$ we have the natural map
$H^*_L(Z,\ome_Z)\to H^*_L(Y,\ome_Y)$ .

We may also consider the
ordinary equivariant homology $H_{*,L}(Y,\CC)$ of $Y$ which is (by definition) equal to
$H^*_{c,L}(Y,\ome_Y)$.
\footnote{Here and in what follows the subscript "c" denotes the appropriate cohomology
with compact support}. In particular we have the natural ("forgetting the supports")
morphism $H_{*,L}(Y,\CC)\to H^*_L(Y,\ome_Y)$.
Also since there is a canonical morphism $\CC_Y\to \ome_Y[-2\dim Y]$
(here $\CC_Y$ denotes the corresponding constant sheaf on $Y$) we get a canonical map
$H^*_L(Y,\CC)\to H^{*-2\dim Y}_L(Y,\ome_Y)$ where $H^*_L(Y,\CC)$
denotes the usual equivariant cohomology of $Y$.

In addition we have the natural pairing
\eq{pairing-raz}
H_{*,L}(Y,\CC)\underset{\calA_L} \ten H^*_L(Y,\CC)\to \calA_L
\end{equation}
(it is clear that all the
(co)homologies in question are modules over $\calA_L$). In particular, we have the integration
functional
$$
\int\limits_Y: H_{*,L}(Y,\CC)\to \calA_L
$$
which is defined as the pairing with $1\in H^*_L(Y,\CC)$. For any closed $L$-invariant subscheme $Z$
of $Y$ the diagram
$$
\begin{CD}
H_{*,L}(Z,\ome_Z)@>>> H_{*,L}(Y,\ome_Y)\\
@V{\int_Z}VV      @VV{\int_Y}V\\
\calA_L@>{\id}>> \calA_L
\end{CD}
$$
is commutative.
%----------------------------------------------------
\ssec{}{Integration on Borel-Moore homology}
We now want to explain that in certain situations one can extend the above integration functional
to the Borel-Moore homology. Assume first that $L$ is a torus and set $Z=Y^L$.
This is a closed subscheme
of $Y$. It is well-known that in this case the natural morphism
$H^*_L(Z,\ome_Z)=H^*(Z,\ome_Z)\ten \calA_L\to
H^*_L(Y,\ome_Y)$ becomes an isomorphism after tensoring with $\calK_L$.

Assume now (until the end of this section) that $Z$ is proper. Then
$H^*_{c,L}(Z,\ome_Z)\simeq H^*_L(Z,\ome_Z)$ we get the map
$H^*_L(Y,\ome_Y)\underset{\calA_L}\ten\calK_L\to H_{*,L}(Y,\CC)\underset{\calA_L}\ten\calK_L$
as the composition
\eq{bm-to-ord}
H^*_L(Y,\ome_Y)\underset{\calA_L}\ten\calK_L\simeq H^*_L(Z,\ome_Z)\underset{\calA_L}\ten\calK_L\\
\simeq H^*_{c,L}(Z,\ome_Z)\underset{\calA_L}\ten\calK_L\to
H_{*,L}(Y,\CC)\underset{\calA_L}\ten\calK_L.
\end{equation}
Setting $H^*_L(Y,\ome_Y)\underset{\calA_L}\ten\calK_L=H^*_{\calK_L}(Y,\ome_Y)$ and
$H^*_L(Y,\CC)\underset{\calA_L}\ten \calK_L=H^*_{\calK_L}(Y,\CC)$
we see that from \refe{pairing-raz} and \refe{bm-to-ord} we get the natural
pairing
\eq{pairing-dva}
H^*_{\calK_L}(Y,\ome_Y)\underset{\calK_L}\ten H^*_{\calK_L}(Y,\CC)\to \calK_L.
\end{equation}
In particular, by pairing with the unit class we get the integration functional
$$
\int\limits_Y:H^*_{\calK_L}(Y,\ome_Y)\to \calK_L.
$$
Since we have the obvious map $H^*_L(Y,\ome_Y)\to H^*_{\calK_L}(Y,\ome_Y)$ the functional
$\int\limits_Y$ makes sense on $H^*_L(Y,\ome_Y)$ too (but of course it still takes values in
$\calK_L$).
This functional is "homogeneous" in the following sense. Assume that we are given
some $\eta\in H^i_L(Y,\ome_Y)$. If $i$ is odd then $\int_Y\eta=0$. If $i$ is even then
$\int_Y\eta$ considered as a rational function on $\grl$ is homogeneous of degree $\frac{i}{2}$.

It is easy to see that all of the above constructions make sense for an arbitrary reductive $L$
provided that the scheme $Y^T$ is proper where $T\subset L$ is a maximal torus.
%------------------------------------------------------------------------------------------------------
\ssec{}{Integration on ordinary cohomology and the integral of the unit class}
By composing $\int_Y$ with the map
$H^*_L(Y,\CC)\to H^{*-2\dim Y}(Y,\ome_Y)$ we may define the integration
functional on $H^*_L(Y,\CC)$ (or $H^*_{\calK_L}(Y,\CC)$). Abusing slightly the notation we shall
denote this functional by the same symbol. This functional also vanishes on the odd part
of $H^*_L(Y,\CC)$ and it sends any class $\psi\in H^{2d}_L(Y,\CC)$ to an element of $\calK_L$
which (considered as a rational function on $\grl$) is homogeneous of degree $d-\dim Y$.
In particular, we may consider
the integral
$$
\int\limits_Y 1_Y
$$
of the unit class $1_Y\in H^*_L(Y,\CC)$. This is a rational function on $\grl$ of degree
$-\dim Y$.

Here is the basic example of the above situation.
Assume that $Y$ is smooth. Let
$\grt$ denote the Lie algebra of $T$. Let also $W_L$ denote  the Weyl group of $L$. Then
$\calA_L=\calA_T^{W_L}$ and $\calK_L=\calK_T^{W_L}$. For every
$y\in Y^T$ let $d_y\in\calA_T$ denote the determinant of the
$\grt$-action on the cotangent space $T_y^*Y$ to $Y$ at $y$. We have
$d_y\neq 0$ if and only if $y$ is an isolated fixed point. Assume
that $Y^T$ is finite. Then it is well known that
    \eq{dh}
    \int\limits_Y 1_Y=\sum\limits_{y\in
    Y^T}\frac{1}{d_y}\in\calK_T^{W_L}=\calK_L.
    \end{equation}

More generally, assume that both $Y$ and $Y^T$ are smooth. Let
$c\in H^*(Y^T,\CC)\ten \calA_T$ denote the equivariant Euler class
of the conormal bundle to $Y^T$ in $Y$. Then it is easy to see that
$c$ is invertible in $H^*(Y^T,\CC)\ten \calK_T$ and we have
$$
\int\limits_{Y}1_Y=\int\limits_{Y^T}\frac{1}{c}.
$$
In general, when $Y$ is singular (and, say, $Y^T$ is finite) the integral $\int\limits_Y 1_Y$
can also be written as a sum over all $y\in Y^T$ of certain "local
contributions" taking values in $\calK_T$. However, in the general
case there is no good formula for these local contributions. One
of the purposes of this paper is to compute $\int\limits_Y 1_Y$ in
certain special cases.

The following lemma follows easily from the definitions and will be used several times in this paper:

\lem{}Let $f:Y_1\to Y_2$ be a proper birational map of $L$-varieties. Assume also that
$\int\limits_{Y_2}1_{Y_2}$ is well-defined. Then $\int\limits_{Y_1}1_{Y_1}$ is also well-defined and we have
$$
\int\limits_{Y_1} 1_{Y_1}=\int\limits_{Y_2}1_{Y_2}.
$$
\elem
%-------------------------------------------------------------------------------------
\ssec{}{The intersection cohomology approach}
We now want to give yet another definition of the above integration based on intersection
cohomology. In fact, in order to do this we will need to make an additional assumption about the
action of $L$ on $Y$ (which will be satisfied for the spaces $^{\bt}\QM^{\theta}_{G,P}$ and
$\calU^{\theta_{\aff}}_{G,P}$ discussed in the introduction).

 Let $Y$ and $L$ be as before and let also $\IC_Y$ denote the intersection cohomology
complex of $Y$. Set
\eq{}
    \IH^*_L(Y)=H^*_L(\IC_Y);\qquad \IH^*_{c,L}=H^*_{c,L}(\IC_Y).
\end{equation}
Both spaces $\IH^*_L(Y)$ and $\IH^*_{c,L}(Y)$ are graded
${\calA_L}$-modules. We have the natural equivariant Poincare
pairing
$$
\la\cdot,\cdot\ra_{Y}:\IH^*_L(Y)\underset{{\calA_L}}\ten
\IH^*_{c,L}(Y)\to{\calA_L}.
$$
This pairing is homogeneous of degree $-\dim Y$. Also we have the
natural (forgetting the supports) morphism
$$
\IH^*_{c,L}(Y)\to \IH^*_L(Y)
$$
of graded ${\calA_L}$-modules.

Set $\IH_{\calK_L}(Y)=\IH_L(Y)\underset{{\calA_L}}\ten{\calK_L}$
and
$\IH_{c,{\calK_L}}(Y)=\IH_{c,L}(Y)\underset{{\calA_L}}\ten{\calK_L}$.
Abusing slightly the notation we shall denote the
corresponding pairing between $\IH_{c,{\calK_L}}(Y)$ and
$\IH_{\calK_L}(Y)$ (which now takes values in ${\calK_L}$) also by
by $\la\cdot,\cdot\ra_{Y}$).

We have the natural morphisms $\CC_Y[\dim Y]\to \IC_Y$ and $\IC_Y\to \ome_Y[-\dim Y]$
which give rise to the maps $H^{*+\dim Y}_L(Y,\CC)\to \IH^*_L(Y)$ and
$\IH_{c,L}^*(Y)\to H^{*-\dim Y}_{c,L}(Y,\ome_Y)=H_{*-\dim Y}(Y,\CC)$
which are compatible with the above pairings in the obvious sense.
In particular we may consider the image of the unit cohomology class
$1_Y$ in $\IH^*_L(Y)$ which we shall denote by the symbol.

Assume now that the map $\IH_{c,\calK_L}\to \IH_{\calK_L}$ is an isomorphism.
Then we may consider the Poincar\'e pairing as a paring
$$
\la\cdot,\cdot\ra_Y:\IH_{\calK_L}(Y)\ten\IH_{\calK_L}(Y)\to\calK_L.
$$

Let
$1_Y\in\IH_{\calK_L}(Y)$ denote the unit cohomology class. For
every $\ome\in\IH_{{\calK_L}}$ we set
$$
\int\limits_Y \ome=\la\ome,1_Y\ra_{Y}\in{\calK_L}.
$$
In particular one may consider the integral
$\int_Y 1$.

Here is a condition which guarantees that the map $\IH_{c,\calK_L}\to
\IH_{\calK_L}$ is an isomorphism.
\lem{}
Let $Y$ be as above. Assume that

a) $Y^T$ is proper;

b) There exists an embedding $\CC^*\subset T$ whose action on $Y$ is
contracting
(this means that the action morphism $\CC^*\x Y\to Y$ extends to a
morphism
$\CC\x Y\to Y$.

Then the natural map $\IH_{c,\calK_L}\to \IH_{\calK_L}$ is an isomorphism.
\elem

In particular, under the assumptions of the lemma one may consider the
integral
$\int_Y 1\in\calK_L$. It is easy to see that this definition of
$\int_Y 1$ coincides with the one discussed above.

%--------------------------------------------------
%----------------------------------------------------------------------
\sec{main}{The main result}
The purpose of this section is to formulate more carefully the results
announced in the introduction.

%----------------------------------------------------------------------------------------------------------------

In this section we are going to work over an arbitrary field $k$ of
characteristic $0$ (later $k$ will become the field $\calK_T$
considered above).

Let $G,\grg$ be as in the introduction - i.e. $\grg$ is a simple Lie
algebra and $G$ is the corresponding simply connected algebraic
group.
The Lie algebra $\grg$ has its
triangular decomposition
$$
\grg=\grn_+\oplus\grh\oplus\grn_-.
$$
Let $T\subset G$ be the corresponding Cartan subgroup. We denote by
$\Lam_G$ the lattice of cocharacters of $T$.
The lattice $\Lam_G$ is spanned by the coroots $\alp:\GG_m\to T$ and
we shall refer to $\Lam_G$ as the coroot lattice. In fact, we let
$\calI_G$ denote the set of verties of the Dynkin diagram of $G$ then
$\Lam_G$ is spanned by the simple coroots $\alp_{\gri}$ for $\gri\in
\calI_G$. We also let
$\Lam_G^{\pos}\subset \Lam_G$ denote the sub-semigroup whose elements
are
linear combinations of positive coroots with non-negative coefficients
(equivalently, linear combinations of $\alp_{\gri}$'s with non-negative coefficients).

Let $P\subset G$ be a parabolic subgroup with Levi component $M$. We
denote by
$\Lam_{G,P}$ the lattice of cocharacters of the group $M/[M,M]$. We
have the natural map $\Lam_G\to \Lam_{G,P}$ and we let
$\Lam_{G,P}^{\pos}$ denote the positive span of the images of
$\alp_{\gri}$'s in $\Lam_{G,P}$ for $\gri\in
\calI_G\backslash\calI_M$.

%---------------------------------------------------------------
\ssec{}{Verma modules and Whittaker vectors: the finite dimensional case}
Let $\chg$ denote the
Langlands
dual algebra of $\grg$. We have the triangular decomposition
$$
\chg=\chn_+\oplus\chh\oplus\chn_-.
$$

Note that the set of vertices of the Dynkin diagram of $\chg$ is equal
to $\calI_G$. We denote by $e_{\gri}\in\chn_+$ and $f_{\gri}\in\chn_-$
(for $i\in\calI_G$)
the corresponding Chevalley generators of $\chg$.
We denote by $\kap:\chg\to\chg$ the Cartan anti-involution sending
$e_{\gri}$ to $f_{\gri}$ and acting on $\chh$ as identity.

Recall that the lattice $\Lam_G$ defined above may
be considered as the coroot lattice of $\chg$. Also we may
think of $\grh$ as the space of all (not necessarily integral) weights
of $\chg$. We let $\rho\in\grh$ denote the half-sum of the positive
roots of $\chg$. Also we denote by $\check{\rho}$ the half-sum of the
positive roots of $\grg$.

Choose now any $\lam\in\grh$.
We denote by $M(\lam)$ the $\chg$-module generated by a
vector $v_{\lam}$ such that:

1) $\chn_-$ annihilates $v_{\lam}$;

2) For any $t\in \chh$ one has $t(v_{\lam})=\lam(t)v_{\lam}$.

3) $v_{\lam}$ is a free generator of $M(\lam)$ over $U(\chn_+)$.

We shall call $M(\lam)$ the {\it Verma module with lowest weight
$\lam$}. We also denote by $\del M(\lam)$ the corresponding dual Verma
module; by the definition it consists of all $\chh$-finite linear
functionals $m:M(\lam)\to k$ with the $\chg$-action defined by
$$
x(m)(v)=m(\kap(x)(v))
$$
for $x\in\chg,v\in M(\lam),m\in \del M(\lam)$.

We have the unique map $\eta_\lam:M(\lam)\to \del M(\lam)$
of $\chg$-modules satisfying $\eta_\lam(v_{\lam})(v_{\lam})=1$.
One can interpet this map as a pairing
$$
\la\cdot,\cdot\ra:M(\lam)\ten M(\lam)\to k
$$
satisfying
$$
\la x(v),w\ra=\la v,\kap(x)(w)\ra
$$
for each $x\in\chg$ and $v,w\in M(\lam)$.
It is well-known that for generic $\lam$ the map $\eta_\lam$ is an
isomorphism.
In fact for this it is enough to require, for example, that
$\la\lam,\alp^{\vee}\ra\not\in\ZZ$
for all coroots $\alp^{\vee}$ of $\chg$.

Let $\hatM(\lam)$ denote the completion of $M(\lam)$ consisting of all
sums
$$
\sum m_\mu \qquad\text{where $m_\mu\im M(\lam)_{\mu}$}
$$
such that for all $\theta\in\Lam_G^{\pos}$ the set
$$
\{ \mu|\ m_\mu\neq 0\quad\text{and}\quad\mu\not\in\lam+\theta+\Lam_G^{\pos}\}
$$
is finite. It is clear that the Lie algebra $\chg$ still acts on
$\hatM(\lam)$.

It is well-known that for any $\lam$ for which the
map $\eta_\lam$ is an isomorphism and for every $\hbar\in k$ such that
$\hbar\neq 0$ there exists unique vector $w(\lam)\in \hatM(\lam), w(\lam)=\sum
w(\lam)_\mu$ (where $w(\lam)_{\mu}$ has weight $\mu$) such that

1) $w(\lam)_{\lam}=v_{\lam}$;

2) For any $\gri\in\calI_G$ we have $f_{\gri}(w)=\frac{1}{\hbar}w$.

We shall refer to $w$ as the {\it Whittaker vector} in $M(\lam)$ (of
course it depends on the choice of $\hbar$; however later this choice
will become canonical).

We also set
$$
w'(\lam)=\sum_{\mu}(-1)^{\la\lam-\mu,\check{\rho}\ra}w(\lam)_{\mu}.
$$
Clearly $w'(\lam)$ satisfies the equation
$f_{\gri}(w'(\lam))=-\frac{1}{\hbar}w'(\lam)$ for all $\gri\in\calI_G$;
$w'(\lam)$ is uniquely characterized by this property if in addition
one requires that the $\lam$-weight part of $w'$ is equal to $v_{\lam}$.
%-----------------------------------------------------------------------------------------------------------------------------
\ssec{}{The affine case}Let us denote by $\grg_{\aff}$ the full
untwisted affine Lie algebra corresponding to $\grg$ (cf. \cite{Kac}).
As a vector space the Lie algebra $\grg_{\aff}$ can be written as
$$
\grg_{\aff}=\grg[x,x^{-1}]\oplus \CC K\oplus \CC d
$$
where $K$ is a central element of $\grg_{\aff}$ and $\ad(d)$ acts on
$\grg[x,x^{-1}]$ as $x\frac{d}{dx}$. The Cartan subalgebra $\grh_{\aff}$ of
$\grg_{\aff}$ is equal to $\grh\oplus \CC K\oplus\CC d$.
We denote by $\calI_G^{\aff}$ the set of vertices of the Dynkin
diagram of $\grg_{\aff}$. Also we denote by $\chg^{\aff}$ the
Langlands
dual Lie algebra of $\grg_{\aff}$. By the definition this is an affine
Lie algebra whose Cartan matrix is the transposed of that of
$\grg_{\aff}$; note that $\chg_{\aff}$ might be a twisted affine Lie algebra.
For each $\gri\in \calI_G^{\aff}$ we denote by
$\alp_{\gri}^{\aff}$ the corresponding simple root of the algebra $\chg_{\aff}$.

The Lie algebras $\grg_{\aff}$ and $\chg_{\aff}$ admit triangular decompositions
$$
\grg_{\aff}=\grn_{+,\aff}\oplus\grh_{\aff}\oplus\grn_{-,\aff}
$$
and
$$
\chg_{\aff}=\chn_{+,\aff}\oplus\chh_{\aff}\oplus\chn_{-,\aff}.
$$
Here $\grh_{\aff}=\grh\oplus\CC d\oplus\CC K$ and
$\chh_{\aff}=\chh\oplus\CC \chd\oplus\CC \chK$.
We have the natural perfect paring between $\grh_{\aff}$ and
$\chh_{\aff}$.

Let $\grh_{\aff}'=\grh\oplus\CC K$ and $\grh''_{\aff}=\grh\oplus\CC d$.
Similarly we define $\chh_{\aff}'$ and $\chh_{\aff}''$. It is well-known that the above perfect pairing
identifies $\grh_{\aff}'$ with the dual space of $\chh_{\aff}''$; it also identifies $\grh_{\aff}''$ with
the dual space of $\chh_{\aff}'$.

We also choose an affine analogue $\rho_{\aff}\in \grh''_{\aff}$ of $\rho$ (cf. \cite{Kac} for a detailed discussion).

Set $\chg_{\aff}'$ to be the subalgebra of $\chg_{\aff}$ consisting of elements whose projection to $\CC \chd$
is equal to $0$. We have
$$
\chg'=\chn_{+,\aff}\oplus\chh_{\aff}'\oplus\grn_{-,\aff}.
$$
For any $\lam_{\aff}=(\lam,\eps)\in \grh_{\aff}''=(\chh_{\aff}')^*$
we let $M(\lam_{\aff})$ denote the corresponding Verma module with lowest weight
$\lam_{\aff}$ over $\chg'_{\aff}$.
In other words this module is freely generated by a vector $v_{\lam_{\aff}}$ subject to the conditions:

1) $t(v_{\lam_{\aff}})=\lam_{\aff}(t)v_{\lam_{\aff}}$ for any $t\in\chh_{\aff}'$.

2) $v_{\lam_{\aff}}$ is annihilated by $\grn_{+,\aff}$.

It is well-known that for generic $\eps$ (more precisely for all but one (critical) value of $\eps$) the
$\chg_{\aff}'$-module structure on $M(\lam_{\aff})$ extends uniquely to an action of the whole $\chg_{\aff}$
such that $d$ annihilates $v_{\lam_{\aff}}$.
Thus for such $\eps$ we may regard $M(\lam_{\aff})$ as a $\chg_{\aff}$-module.

The Lie algebra $\chg_{\aff}$ carries canonical anti-involution $\kap_{\aff}$ (interchanging
$\chn_{+,\aff}$ and $\chn_{-,\aff}$ and acting as identity on $\chh_{\aff}$). Similarly to the affine case there is
a unique pairing
$$
\la\cdot,\cdot\ra: M(\lam_{\aff})\ten M(\lam_{\aff})\to \CC
$$
satisfying
$$
\la x(v),w\ra=\la v,\kap_{\aff}(x)(w)\ra
$$
for each $x\in\chg$ and $v,w\in M(\lam_{\aff})$ and such that
$$
\la v_{\lam},v_{\lam}\ra=1.
$$
Of course all the above constructions make sense when the base field $\CC$ is replaced by any field $k$ of
characteristic 0 (later $k$ will become certain field of rational functions).

As before one may define a completion $\hatM(\lam_{\aff})$ of $M(\lam_{\aff})$. Choose $\hbar\in k^*$.
Then for generic
$\lam_{\aff}$ this completion will contain unique vector $w(\lam_{\aff})=\sum w(\lam_{\aff})_{\mu}$
(here $\mu\in \grh_{\aff}''$) subject to the conditions

1) $w(\lam_{\aff})_{\lam_{\aff}}=v_{\lam_{\aff}}$;

2) For any $\gri\in \calI_{G}^{\aff}$ we have $f_{\gri}(w(\lam_{\aff}))=\frac{1}{\hbar}w(\lam_{\aff})$.

Again we shall call $w(\lam_{\aff})$ the {\it Whittaker vector} of
$M(\lam_{\aff})$.
Similarly one can define the vector $w'(\lambda_{\text{aff}}$.

%------------------------------------------------------------------------------------------------------
\ssec{}{The main result}
Set now
$$
\IH_{G,P}=\IH_{\calK_{M\x\CC^*}}(\QM^\theta_{G,P})
\quad\text{and}\quad
\IH_{G,P}^{\theta_{\aff}}=\IH_{\calK_{M\x (\CC^*)^2}}(\calU^{\theta_{\aff}}_{G,P}).
$$
Also we denote
$$
\IH_{G,P}=\bigoplus\limits_{\theta\in
\Lam_{G,P}^{\pos}}\IH_{G,P}^{\theta}
\quad\text{and}\quad
\IH_{G,P}^{\aff}=\bigoplus\limits_{\theta_{\aff}\in
\Lam_{G,P}^{\aff,\pos}}\IH_{G,P}^{\theta_{\aff}}.
$$

We shall further restrict ourselves to the case $P=B$ (where $B$ as before denotes a Borel subgroup of $G$).
In this case $\IH_{G,B}$ is a vector space over the field $\calK_{T\x \CC^*}$ which we may identify with the field
of rational functions on $\grh\x\CC$. We shall denote a typical element in this product by $(a,\hbar)$.
Similarly $\IH_{G,B}^{\aff}$ is a vector space over the field $\calK_{T\x (\CC^*)^2}$ which we may identify with the
field of rational functions on $\grh\x\CC^2$. We shall typically denote an element of the latter product
as $(a,\eps,\hbar)$. Note that we may think of $(a,\eps)$ as an element of $\grh_{\aff}'$. Note also that we have
a canonical element $\hbar$ in both $\calK_{T\x\CC^*}$ and $\calK_{T\x (\CC^*)^2}$.

We denote by $\la \cdot,\cdot\ra_{G,B}^{\theta}$ (resp. by $\la
\cdot,\cdot\ra_{G,B}^{\theta_{\aff}}$) the Poincare pairing on
$\IH_{G,B}^{\theta}$ (resp. on $\IH_{G,B}^{\theta_{\aff}}$). As in the
introduction we let $\la\cdot,\cdot\ra_{G,B}$ be the direct sum of the
pairings
$(-1)^{\la \theta,\check{\rho}\ra}\la
\cdot,\cdot\ra_{G,B}^{\theta}$. Similarly we let
$\la\cdot,\cdot\ra_{G,B}^{\aff}$ denote the direct sum of the pairings
$(-1)^{\la \theta_{\aff},\check{\rho}_{\aff}\ra}\la \cdot,\cdot\ra_{G,B}^{\theta_{\aff}}$.
%-------------------------------------------------------------------------------------------------
\th{main} There exists a
natural action of $\chg$ (resp. of $\chg_{\aff}$) on $\IH_{G,B}$
(resp. on $\IH_{G,B}^{\aff}$) such that:
\begin{enumerate}
\item
$\IH_{G,B}^{\theta}$
has weight $\frac{a}{\hbar}+\rho+\theta$ with respect to the action of $\chg$.
Similarly,
$\IH_{G,B}^{\theta_{\aff},\aff}$
has weight $\frac{(a,\eps)}{\hbar}+\rho_{\aff}+\theta_{\aff}$ with respect to
the action of $\chg_{\aff}$
\footnote{here we regard $(a,\eps)$ as an element of $\grh_{\aff}$ by means of the embedding
$\grh_{\aff}'\subset\grh_{\aff}$}.
%%%%%%%%%%%%%%%%%
\item
The pairing $\la\cdot,\cdot\ra_{G,B}$ is $(\chg,\kap)$-invariant.
Similarly, the pairing $\la\cdot,\cdot\ra_{G,B}^{\aff}$ is
$(\chg_{\aff},\kap_{\aff})$-invariant.
%%%%%%%%%%%%%%%%%%%%%%%%%%%%%%%%%%%%%%%%%%%%
\item The $\chg$-module $\IH_{G,B}$ is
isomorphic to $M(\frac{a}{\hbar}+\rho)$ as a $\chg$-module. Similarly,
$\IH_{G,B}^{\aff}$ is isomorphic to
$M_{\aff}(\frac{(a,\eps)}{\hbar}+\rho_{\aff})$.

\item
Under the above identifications the vector $\sum
1_{\QM_{G,B}^\theta}$ (resp. $\sum
1_{\calU_{G,B}^{\theta_{\aff}}}$) is the Whittaker vector
$w(\frac{a}{\hbar}+\rho)$
(resp. $w(frac{(a,\eps)}{\hbar}+\rho_{\aff})$) in
$\hatM(\frac{a}{\hbar}+\rho)$ (resp. in
$\hatM(\frac{(a,\eps)}{\hbar}+\rho_{\aff})$).
\end{enumerate}
\eth
In what follows we let $a_{\aff}$ denote the pair $(a,\eps)$ (in the affine case).
%-----------------------------------------------------------------------------------------------
\cor{main-c}
We have
$$
\calZ_{G,B}=\grq^{-\frac{a}{\hbar}-\rho}\la w(\frac{a}{\hbar}+\rho),\grq\cdot w'(\frac{a}{\hbar}+\rho)\ra
\footnote{Of course the torus $\chT$ does not act on $M(\frac{a}{\hbar})$ (only its Lie algebra acts).
Thus in order to make sense of this formula we must (as before) write $q=e^t$ where $t\in \chh$ and thus
the equality becomes an equality of formal Fourier series in $t$.}
$$
Similarly,
$$
\calZ_{G,B}^{\aff}=\grq^{-\frac{a}{\hbar}-\rho} \la w(\frac{a_{\aff}}{\hbar}+\rho_{\aff}),
\grq_{\aff}\cdot w'(\frac{a_{\aff}}{\hbar}+\rho_{\aff})\ra.
$$
\ecor

\ssec{}{The quadratic Toda Hamiltonian}
Fix now an invariant quadratic form $(\cdot,\cdot)$ on $\grg$ so that long roots on
$\chg$ have squared norm 2. Denote by $\Del$ the corresponding
Laplacian on $\chh$. As before we set $\grq=e^t$ where $t\in \chh$ and
$\grq_{\aff}=e^{t_{\aff}}$. Thus we may write $t_{\aff}=(t,\ln Q)$.
\cor{toda}
\begin{enumerate}
\item
The function $e^{\frac{a(t)}{\hbar}}\calZ_{G,B}$ is an eigen-function of the
quantum Toda Hamiltonians with the "eigen-value" defined by $a$
(cf. \cite{etingof} and references therein for the definition).
In particular,  the function $\calZ_{G,B}$ satisfies the equation
$$
(\hbar^2\Del-2\sum\limits_{\gri\in\calI_G}e^{\alp_{\gri}(t)})
(e^{\frac{a(t)}{\hbar}}\calZ_{G,B})=(a,a)e^{\frac{a(t)}{\hbar}}\calZ_{G,B}
$$
\item
The function $\calZ_{G,B}^{\aff}$  satisfies
\eq{non-stat}
(\eps Q\frac{\partial}{\partial
Q}+\hbar^2\Del-2\sum\limits_{\gri\in\calI_G^{\aff}}e^{\alp^{\aff}_{\gri}(t_{\aff})})
e^{\frac{a(t)}{\hbar}}\calZ_{G,B}^{\aff}=(a,a)e^{\frac{a(t)}{\hbar}}\calZ_{G,B}.
\end{equation}
\end{enumerate}
\ecor
The operator in the left hand side of \refe{non-stat} may be called a
{\it non-stationary analogue} of the quadratic affine Toda hamiltonian
(the latter will be obtained from it by removing the $\eps Q\frac{\partial}{\partial
Q}$-term).
%---------------------------------------------------------------------------------------------------
\sec{sl2}{Example: $G=SL(2)$} In this section we would like to prove
\refc{main-c} explicitly for $G=SL(2)$ (in the non-affine case).
So, we shall assume that $G=SL(2)$ until the end of this section.
%---------------------------------------------------------------------------------------------------------
\ssec{}{The space $^{\bt}\calM_{G,B}^d$}Since $G=SL(2)$ we may
identify $\Lam_G$ with $\ZZ$ and we shall typically denote an element
of $\Lam_G$ by $d$ (instead of $\theta$).
Also in this case $\calG_{G,B}=\PP^1$ (we choose this identification
in such a way that $e_{G,B}$ corresponds to $\infty\in\PP^1$) and thus $^{\bt}\calM_{G,B}^d$
is just the space of maps $f:\PP^1\to\PP^1$ of degree $d$ which send
$\infty$ to $\infty$. This is the same as the set of pairs $(g,h)$ of
polynomials
in one variable such that:

1)$g$ and $h$ have no common roots

2)$g(z)=z^d+\text{lower order terms}$

3)$\deg h< d$.
%-------------------------------------------------------------------------------------------
\ssec{}{The space $^{\bt}\QM_{G,B}^d$ and the integral of $1$}
The ``zastava'' space $^{\bt}\QM^d_{G,B}$
consists of all pairs $g,h$ as above but with condition 1 removed
(in particular $h$ is allowed to be identically $0$). Thus
$^{\bt}\QM_{G,B}^d$ is isomorphic to the affine space $\AA^{2d}$. It inherits the
natural action of $T\x \CC^*$ with unique fixed point
corresponding to $g=z^d,h=0$.

Now for each $d$ we set
$$
A_d:=\int\limits_{^{\bt}\QM_{G,B}^d} 1\quad \in\calK_{T\x\CC^*}
$$
(the integral is taken in the $T\x\CC^*$-equvariant cohomology). Note that $T=\CC^*$ in our case; thus we shall
identify $\calK_{T\x\CC^*}$ with the field of rational functions in the variables $a,\hbar$.
Since $\calM^d$ is smooth and we have just one fixed point it is
easy to compute this integral explicitly by using \refe{dh} and we get
$$
A_d=\frac{1}{d!\hbar^d\prod\limits_{i=1}^d(a+i\hbar)}.
$$

%------------------------------------------------------------------------------
\ssec{}{Verma modules and Whittaker vectors}
Let us denote by $e,h,f$ the standard basis for $sl(2)\simeq \chg$. For
every scalar $\lam$ we denote by $M(\lam)$ (as before) the Verma module over
$sl(2)$ with {\it lowest weight} $\lam$.
Explicitly this means that $M(\lam+1)$ has a basis $\{ m_d\}_{d\geq
0}$ in which the action is written in the following way:
$$
h(m_d)=(\lam+2d+1)m_d;\quad e(m_d)=m_{d+1};\quad
f(m_d)=-d(\lam+d)m_{d-1}
$$
(thus $f$ acts nilpotently and $\CC[e]$ acts freely).

The anti-involution $\kap:sl(2)\to sl(2)$ is written in the following way:
$$
\kap(e)=f;\quad \kap(f)=e\quad \kap(h)=h.
$$

Recall that $\la\cdot,\cdot\ra$ denotes the unique $(sl(2),\kap)$-invariant form on $M(\lam+1)$
satisfying $\la m_0, m_0\ra=1$.

Explicitly this form is
written in the following way:
$$
\langle m_d,m_{d'}\rangle=\del_{d,d'} (-1)^d
d!\prod\limits_{i=1}(\lam+i).
$$
We now consider the case $\lam=\frac{a}{\hbar}$. In this case define
$$
w_d=\frac{v_d}{d!\prod\limits_{i=1}^d (a+i\hbar)}.
$$
Thus
$$
f(w_d)=-\frac{1}{\hbar}w_{d-1}.
$$
Set also $w=\sum (-1)^d w_d$, $w'=\sum w_d$. Thus $f(w)=\frac{1}{\hbar}w$ and
$f(w')=-\frac{1}{\hbar}w'$ - i.e. $w=w(\frac{a}{\hbar}+1)$ and $w'=w'(\frac{a}{\hbar}+1)$.

On the other hand we have
\eq{nahui}
\la w_d,w_d\ra =\frac{\la v_d,v_d\ra}{(d!\prod\limits_{i=1}^d (a+i\hbar))^2}=\\
\frac{(-1)^d d!\prod\limits_{i=1}^d(\frac{a}{\hbar}+i)}
{(d!\prod\limits_{i=1}^d (a+i\hbar))^2}=
\frac{(-1)^d}{d!\hbar^d\prod_{i=1}^d(a+i\hbar)}=(-1)^d A_d.
\end{equation}
This is equivalent to \refc{main-c} in this case.
%-----------------------------------------------------------------------------------------------------------
\ssec{}{The geometric interpretation of the action}
To complete the picture for $G=SL(2)$ let us explain the geometric
origin of the operators $(e,h,f)$ on $\IH_{G,B}$. Since $e$ is adjoint to $f$ with respect to $\la\cdot,\cdot\ra$
and $[e,f]=h$ it follows that it is enough to give a construction of $f$.

For each $d\geq 0$ let $Y_d\subset  {^{\bt}\QM_{G,B}^d}\x {^{\bt}\QM_{G,B}^{d+1}}$
consisting of pairs $((g_1,h_1),(g_2,h_2))$  subject to the following
condition:

there exists $z_0\in \CC$ so that $g_2(z)=(z-z_0)g_1(z)$ and $h_2(z)=(z-z_0)h_1(z)$.

We denote by $\pr_1$ the projection from $Y_d$ to
${^{\bt}\QM_{G,B}^d}$ and by $\pr_2$ the corresponding projection to
${^{\bt}\QM_{G,B}^{d+1}}$.

The variety $Y_d$ viewed as a correspondence from
${^{\bt}\QM_{G,B}^{d+1}}$ to ${^{\bt}\QM_{G,B}^d}$ defines a
$\calK_{T\x \CC^*}$-linear map
$H^*_{\calK_{T\x\CC^*}}({^{\bt}\QM_{G,B}^{d+1}})\to H^*_{\calK_{T\x
\CC^*}}({^{\bt}\QM_{G,B}^d})$.
Since we are ignoring the grading and since in this case the varieties
${^{\bt}\QM_{G,B}^d}$
are smooth we have $H^*_{\calK_{T\x
\CC^*}}({^{\bt}\QM_{G,B}^d})=\IH^*_{\calK_{T\x \CC^*}}({^{\bt}\QM_{G,B}^d})$
and thus we let
$f:\IH^*_{\calK_{T\x \CC^*}}({^{\bt}\QM_{G,B}^{d+1}})\to
\IH^*_{\calK_{T\x \CC^*}}({^{\bt}\QM_{G,B}^d})$ be the operator defined
by $Y_d$ as a correspondence. Since  $Y_d$ is clearly isomorphic to
${^{\bt}\QM^d_{G,B}}\x\AA^1$  to see that $f(1_{G,B}^{d+1})=\frac{1^d_{G,B}}{\hbar}$.

We now define $h:\IH^*_{\calK_{T\x \CC^*}}({^{\bt}\QM_{G,B}^{d}})\to
\IH^*_{\calK_{T\x \CC^*}}({^{\bt}\QM_{G,B}^d})$ to be the operator of
multiplication
by $\frac{a}{\hbar}+2d+1$ and $e:\IH^*_{\calK_{T\x \CC^*}}{^{\bt}(\QM_{G,B}^{d}})\to
\IH^*_{\calK_{T\x \CC^*}}({^{\bt}\QM_{G,B}^d})$ to be the conjugate of
$f$ with respect to the pairing $\la\cdot,\cdot\ra_{G,B}$.

We now have to check that the operators $e,h,f$ defined in the above
way satisfy the relations of the Lie algebra $sl(2)$. The relations
$[h,e]=2e$ and $[h,f]=-2f$ are obvious. Thus we have to check only
that $[e,f]=h$. This is equivalent to \refe{nahui}
(since $A_d=(-1)^d \la 1_{G,B}^d, 1_{G,B}^d\ra_{G,B}$).

%%%%%%%%%%%%%%%%%%%%%%%%%%%%%%%%%%%%%%%%%%%%%%%%%%%%%%%%%%%%%%%%%%%%%%%%%%%%%%%%%%%%%%%%%%%%%%%%%%%%%%%%%%%%%%%%%%
%%%%%%%%%%%%%%%%%%%%%%%%%%%%%%%%%%%%%%%%%%%%%%%%%%%%%%%%%%%%%%%%%%%%%%%%%%%%%%%%%%%%%%%%%%%%%%%%%%%%%%%%%%%%%%%%%%
\sec{action}{Construction of the action}
%----------------------------------------------------------------------------------
In this section we sketch the proof of \reft{main} assuming that the reader is familiar with the contents
of \cite{ffkm}.

\ssec{}{Plan of the proof}
First we are going to discuss the finite case.

Recall that $\chg$ has the triangular decomposition
$$
\chg=\chn_+\oplus\chh\oplus\chn_-.
$$
The construction of the $\chg$-action on $\IH_{G,B}$ will consist of the following 3 steps:

$\bullet$ Step 1: Construct an action of $\chn_+$ on $\IH_{G,B}$.

$\bullet$ Step 2: Use the  pairing
$\la\cdot,\cdot\ra_{G,B}:\IH_{G,B}\underset{\calK_{T\x\CC^*}}\ten\IH_{G,B}\to\calK_{T\x\CC^*}$
to construct the action of $\chn_-$.

$\bullet$ Step 3: Check that the two actions of $\chn_+$ and of $\chn_-$ generate an action of $\chg$ and that the
resulting $\chg$-module is isomorphic to $M(\frac{a}{\hbar}+\rho)$.

$\bullet$ Step 4: Check that for every Chevalley generator $f_{\gri}$ corresponding to some
$\gri\in\calI_G$ and every $\theta\in\Lam_G^{\pos}$ we have
\eq{step4}
f_{\gri}(1^{\theta}_{G,B})=\frac{1}{\hbar}1^{\theta-\alp_{\gri}}_{G,B}.
\end{equation}
In fact steps 1-3 come from almost a word-by-word repetition of the contents of Section 4 of \cite{ffkm}
and step 4 is equivalent to the computations of the previous section.
%----------------------------------------------------------------------------------
\ssec{}{Step 1}Recall the definition of the {\it local Ext-algebra} from \cite{ffkm}.
For every $\theta\in\Lam_G^{\pos}$ we denote (following \cite{ffkm}) by $\AA^{\theta}$
the corresponding partially symmetrized power of $\AA^1$ (i.e. elements of $A^{\theta}$ are
formal sums $\sum \theta_i x_i$ with $\theta_i\in\Lam_G^{\pos}$, $x_i\in \AA^1$ and
$\sum\theta_i =\theta$). In \cite{ffkm} the authors construct a closed embedding
$\AA^{\theta}\hookrightarrow {^{\bt}\QM^{\theta}_{G,B}}$ whose image is denoted
by $\partial^{\theta}({^{\bt}\QM^{\theta}_{G,B}})$ (in fact this is the closed stratum of a certain
canonical stratification of ${^{\bt}\QM^{\theta}_{G,B}}$).

Now we set (cf. Section 2 of \cite{ffkm})
$$
\calA_{\theta}=
\Ext ^*_{{^{\bt}\QM^{\theta}_{G,B}}}
(\IC_{\partial^{\theta}({^{\bt}\QM^{\theta}_{G,B}})},\IC_{{^{\bt}\QM^{\theta}_{G,B}}})
$$
and define
$$
\calA=\bigoplus\limits_{\theta\in\Lam_G^{\pos}}\calA_{\theta}.
$$
In \cite{ffkm} the authors define in a geometric way a Hopf algebra structure on
$\calA$ and show that it is isomorphic (as a Hopf algebra) to $U(\chn_+)$.

In section 4 of \cite{ffkm} the authors define an action of $\calA$ on
$\oplus_{\theta\in\Lam_G^{\pos}}\IH^*(\QM_{G,B}^{\theta})$. A word-by-word repetition of that
construction defines an action of $\calA$ on $\IH_{G,B}$. Thus we get an action of the Lie
algebra $\chn_+$ on $\IH_{G,B}$.

%-----------------------------------------------------------------------------------------
\ssec{}{Steps 2,3 and 4}We already have the action of the operators $f_{\gri}$ on $\IH_{G,B}$.
We define the action of every $t\in \chh$ on $\IH_{G,B}$ by requiring that
$t$ acts on $\IH_{G,B}^{\theta}$ by means of multiplication by
$\la t,\frac{a}{\hbar}+\rho+\theta\ra$ on $\IH_{G,B}^{\theta}$.
Also we define $e_{\gri}$ to be the adjoint operator to $f_{\gri}$ with respect to
$\la\cdot,\cdot\ra_{G,B}$. We need to check the following relations:
$$
[e_{\gri},f_{\gri}]=h_{\gri};\quad [e_{\gri},f_{\grj}]=0\quad\text{for $\gri\neq\grj$}
$$
The first relation is proved exactly as the Proposition in section 4.9 of \cite{ffkm} and
the second one is proved exactly as Proposition 4.8 of \cite{ffkm}.

In order to check that $\IH_{G,B}$ is isomorphic to $M(\frac{a}{\hbar}+\rho)$ as a $\chg$-module
we may argue (for example) as follows. It follows from the main theorem of section 3 of \cite{ffkm}
(cf. section 3.1 of \cite{ffkm}) that $\IH_{G,B}$ and $M(\frac{a}{\hbar}+\rho)$ have the same character.
Also, it is clear that $\IH_{G,B}$ lies in the category $\calO$ for
$\chg$ (in other words, $\chh$ acts on it in a semi-simple way with
finite-dimensional weight spaces and $\grn_-$ acts locally nilpotetnly). However, it is easy to see
that $M(\frac{a}{\hbar}+\rho)$ is the only module in this category with such character (in fact, it is
the unique module in the category $\calO$ with the lowest weight $\frac{a}{\hbar}+\rho$).

Finally, we need to check \refe{step4}. However, using the same argument as in 4.9 of \cite{ffkm} we can
reduce it to the case $G=SL(2)$ and thus it becomes equivalent to the calculation of the previous
section.
%-------------------------------------------------------------------------------------
\ssec{}{The affine case}Let us explain how to extend the above argument to the affine case.
Here we need to use the results of \cite{bfg} instead of \cite{ffkm}. First of all the algebra
$\calA$ is defined exactly as above using the varieties $\calU_{G,B}^{\theta_{\aff}}$ instead
of $^{\bt}\QM_{G,B}^{\theta}$ (one has to use the stratification of $\calU_{G,B}^{\theta_{\aff}}$
defined in section 10 of \cite{bfg}). The definition of the Hopf algebra structure is a word-by-word
repetition of the construction in section 2 of \cite{ffkm}. The proof of the isomorphism
$\calA\simeq U(\chn_{+,\aff})$ is the same as in section 4 of \cite{ffkm} except that one has to use
the results of section 16 of \cite{bfg} for counting the dimension of $\calA_{\theta}$
 instead of section 3 of \cite{ffkm}. The rest is proved exactly as before (except that
 in order to prove the isomorphism $\IH_{G,B}^{\aff}\simeq M(\frac{a_{\aff}}{\hbar}+\rho_{\aff})$
 one has to use again section 16 of \cite{bfg} instead of section 3 of \cite{ffkm} while computing the character
 of $\IH_{G,B}^{\aff}$).

%----------------------------------------------------------------------------------
\sec{quantum cohomology}{Proof of \reft{j versus z}}
%-------------------------------------------------------------------------------------
The contents of this section are mostly due to A.~Givental. We include it for the sake of
completeness.

\ssec{}{Some reformulations}
In this section we prove \reft{j versus z}. The arguments here are very close to
those in \cite{giv-lee}. Fix $G,P$ as in the formulation of the
theorem. We need to prove that
for each $\beta$ we have
$$
\int\limits_{^{\bt}\QM^{\beta}_{G,P}}1_{G,P}^{\beta}=(\ev_\beta)_*(\frac{1}{\hbar(c_{\beta}+\hbar)}
$$
under the conventions explained before the formulation of \reft{j
versus z}. Let $^{\bt}\ocM^{\beta}_{0,1}(G,P)$ denote the
space of based stable maps from a genus 0 curve with one marked point to $\calG_{G,P}$;
in other words
$^{\bt}\ocM^{\beta}_{0,1}(G,P)$
is a closed sub-scheme in
$\ocM^{\beta}_{0,1}(G,P):=\ocM^{\beta}_{0,1}(\calG_{G,P})$ corresponding to maps sending the
marked point to $e_{G,P}$.

It is now clear that the restriction map gives rise to an isomorphism
$H^*_G(\ocM^{\beta}_{0,1}(G,P),\CC)\simeq
H^*_M(^{\bt}  \ocM^{\beta}_{0,1}(G,P),\CC)$.
Moreover, the following diagram commutes:
$$
\begin{CD}
H^*_G(\ocM^{\beta}_{0,1}(G,P),\CC)@>>>
H^*_M(^{\bt}\ocM^{\beta}_{0,1}(G,P),\CC)\\
@V{(\ev_{\beta})_*}VV
@VV{\int\limits_{^{\bt}\ocM^{\beta}_{0,1}(G,P)}}V\\
H_G^*(\calG_{G,P},\CC)@>>> H_M^*(pt,\CC)
\end{CD}
$$
where the last map is just the restriction to $e_{G,P}$.

In particular we see that under the above identifications we have
$$
(\ev_\beta)_*(\frac{1}{\hbar(c_{\beta}+\hbar)})=
\int\limits_{^{\bt}\ocM(G,P)}\frac{1}{\hbar({c_{\beta}}|_{^{\bt}\ocM(G,P)}+\hbar)}.
$$
Our next goal will be to rewrite the last integral as the integral of
1 over certain space.
%--------------------------------------------------------------------------------------------------
\ssec{}{Graph spaces}
For every $\theta\in\Lam_{G,P}^{\pos}$ let us introduce
the {\it graph} spaces
$$
\Gra^{\beta}_{G,P}:=\ocM_{0,0}^{\beta,1}(\calG_{G,P}\x \PP^1).
$$
A point of $\Gra^{\beta}_{G,P}$ is a stable map $C\to \calG_{G,P}\x
\PP^1$
which has degree $\theta$ over the first multiple and degree one
over the second one. In the case when $C$ is non-singular the second
projection gives an isomorphism between $C$ and $\PP^1$ and thus our
map is actually equal to the graph of some map $\PP^1\to \calG_{G,P}$
of
degree $\beta$. Thus $\calM^\beta_{G,P}$ is an open subset of
$\Gra^{\beta}_{G,P}$.

We denote by $^{\bt}\Gra^{\beta}_{G,P}$ the corresponding ``based''
version of the graph space; it consists of all stable maps
$C\to \calG_{G,P}$ such that

(i) over some neighbourhood of $\infty\in\PP^1$ the projection $C\to
\PP^1$
is an isomorphism. In this case we denote by $\infty_C$ the (unique) preimage
of $\infty$ in $C$.

(ii) The image of $\infty_C$ in $\calG_{G,P}\x \PP^1$ is $(e_{G,P},\infty)$.

Note that $^{\bt}\Gra^{\beta}_{G,P}$ is not a closed subset of
$\Gra^{\beta}_{G,P}$ (it is only locally closed).
%------------------------------------------------------------------------
\ssec{}{The $\CC^*$-action}
The group $\CC^*$ acts naturally on $\Gra^{\beta}_{G,P}$ (via its
action on $\PP^1$). The fixed points of this action are discussed in \cite{giv-lee};
it is shown in {\it loc. cit.} that we have
$$
(\Gra^{\beta}_{G,P})^{\CC^*}=\bigcup\limits_{\beta_1+\beta_2=\beta}
{\ocM^{\beta_1}_{0,1}(G,P)}\x
{\ocM^{\beta_2}_{0,1}(G,P)}
$$
This isomorphism is defined as follows. Assume that we have a stable map
$f:C\to \calG_{G,P}\x\PP^1$ of degree $(\beta,1)$. Then there exists unique
component $C'$ of $C$ which maps isomorphically to the second multiple. Assume now
that $f$ is $\CC^*$-invariant. Then the composition $p_1\circ f|_{C'}$ must be constant
(here $p_1$ is the projection $\calG_{G,P}\x\PP^1\to\calG_{G,P}$).
Let $C''$ be the union of all the irreducible components of $C$ which are different from $C'$.
Then the composition map $C''\to\calG_{G,P}\x\PP^1\to \PP^1$ sends every connected component
of $C''$ either to $0$ or to $\infty$. We set $C''_0=(p_2\circ f|_{C''})^{-1}(0)$ and
$C''_{\infty}=(p_2\circ f|_{C''})^{-1}(\infty)$. We also set $c_0=C''_0\cap C'$ and
$c_{\infty}=C''_{\infty}\cap C'$.
Thus one can write $C''=C''_0\cup C''_{\infty}$ and the restriction of
$p_1\circ f$ to either $(C''_0,c_0)$ or $(C''_{\infty},c_{\infty})$ is a stable map.
In particular, if $(C,f)\in {} ^{\bt}\Gra^{\beta}_{G,P}$ then $C''=C''_{0}$
(since the unique point of $C'$ lying over $\infty\in\PP^1$ must be a non-singular
point of $C$) and the map $(C,f)\mapsto (C'',c_0,p_1\circ f|_{C''})$ defines an isomorphism
$$
(^{\bt}\Gra^{\beta}_{G,P})^{\CC^*}\simeq {^{\bt}\ocM^{\beta}_{0,1}(G,P)}
$$
Moreover, it is easy to see (cf. \cite{giv-lee}, proof of Theorem 2 or \cite{giv}) that the conormal bundle to
$^{\bt}{\ocM^{\beta}_{0,1}(G,P)}$ in $^{\bt}\Gra^{\beta}_{G,P}$ is of rank 2 and its
$\CC^*$-equivariant first Chern class is
$\hbar(c_{\beta}+\hbar)$. Therefore we have
%-----------------------------------------------------------------------------------
\eq{hhh}
\int\limits_{^{\bt}\Gra^{\beta}_{G,P}}
1=\int\limits_{^{\bt}\ocM^{\beta}_{0,1}(G,P)}
\frac{1}{\hbar(c_{\beta}+\hbar)}
\end{equation}
%------------------------------------------------------------------------------------
\ssec{}{End of the proof}It follows from \refe{hhh} that in order
to prove \reft{j versus z} it is enough to show that
\eq{}
\int\limits_{^{\bt}\Gra^{\beta}_{G,P}}1=\int\limits_{^{\bt}\QM^{\beta}_{G,P}}1
\end{equation}
For this it is enough to show the existence of a proper birational $M$-equivariant map
$^{\bt}\Gra^{\beta}_{G,P}\to {^{\bt}\QM^{\beta}_{G,P}}$. This is proved in \cite{ffkm}.

\end{document}